\documentclass[a4paper,12pt]{article}
\usepackage[a4paper,margin=1in,head=14.5pt]{geometry}
\usepackage[utf8]{inputenc}
\usepackage[T1]{fontenc}
\usepackage[english]{babel}
\usepackage[hidelinks,hypertexnames=false]{hyperref}
\usepackage{graphicx}
\usepackage{amsmath}
\usepackage{amssymb}
\usepackage{amsthm}
\usepackage{cleveref}
\usepackage{thmtools}
\usepackage{mathtools}
\usepackage{bbm}
\usepackage{enumerate}
\usepackage{enumitem}
\usepackage{color}
\numberwithin{equation}{section}

\usepackage{scrlayer-scrpage}
\crefname{subsection}{subsection}{subsections}
\crefname{subsubsection}{subsubsection}{subsubsections}
\crefname{figure}{figure}{figures}

\let\originalleft\left
\let\originalright\right
\renewcommand{\left}{\mathopen{}\mathclose\bgroup\originalleft}
\renewcommand{\right}{\aftergroup\egroup\originalright}

\makeatletter
\newcommand{\myitem}[1]{%
\item[#1]\protected@edef\@currentlabel{#1}%
}
\makeatother

\newcommand{\R}{{\mathbb R}}
\newcommand{\N}{{\mathbb N}}

\newcommand{\C}{{\mathbb C}} 
\newcommand{\Cstar}{\C^{*}}

\DeclareMathOperator{\Div}{div}

\DeclareMathOperator{\id}{id}

\newcommand{\opA}{\mathbb A}
\newcommand{\opl}[1][s]{\ell(s)}
\newcommand{\opL}[1][s]{\mathcal L(#1)}
\newcommand{\opN}[1][s]{\mathcal N(#1)}
\newcommand{\opLdual}[1][s]{\mathcal L(#1)}
\newcommand{\opNdual}[1][s]{\mathcal N(#1)}
\newcommand{\opNddual}[1][s]{\mathcal N_{\mathrm{ext}}({#1})}
\newcommand{\opS}[1][s]{\operatorname{\mathcal S}(#1)}
\newcommand{\opD}[1][s]{\operatorname{\mathcal D}(#1)}
\newcommand{\opSD}[1][s]{\operatorname{\mathcal{E}}_{\mathrm{D,ext}}(#1)}
	\newcommand{\opSN}[1][s]{\operatorname{\mathcal{E}}_{\mathrm{N,ext}}(#1)}
\newcommand{\opSj}{\operatorname{\mathcal S}_j}
\newcommand{\opDj}{\operatorname{\mathcal D}_j}
\newcommand{\DtN}[1][s]{\operatorname{DtN}(#1)}
\newcommand{\DtNLaplace}[1][s]{\operatorname{DtN}^{\Delta}}
\newcommand{\DtNLaplaceLog}[1][s]{\operatorname{DtN}^{\Delta}_{\mathrm{log}}}
\newcommand{\NtD}[1][s]{\operatorname{NtD}(#1)}

\newcommand{\opV}[1][s]{\mathsf{V}(#1)}
\newcommand{\opK}[1][s]{\mathsf{K}(#1)}
\newcommand{\opKdual}[1][s]{\mathsf{K}'(#1)}
\newcommand{\opW}[1][s]{\mathsf{W}(#1)}

\newcommand{\opAvar}[1][s]{\mathsf{A}({#1})}

\newcommand{\opC}[1][s]{\mathsf{C}(#1)}
\newcommand{\opCj}[1][s]{\mathsf{C}_j(#1)}

\newcommand{\HGammaD}{H^1_{\Gamma_{\mathrm D}}}
\newcommand{\HdualGammaD}{\Hdual_{\Gamma_{\mathrm D}}}
\newcommand{\IMatrix}{\mathbb I}
\newcommand{\IGebiet}{G}
\newcommand{\IGebietR}{G_R}
\newcommand{\IGebietInt}{G_0}

\newcommand{\Interface}{\Gamma}
\renewcommand{\dim}{{n}}
\newcommand{\conj}[1]{\overline{#1}}
\newcommand{\Hdual}{\widetilde{H}^{-1}}
\newcommand{\Xmult}{\mathbb X}

\newcommand{\nomathfrak}{}

\renewcommand{\Re}{\operatorname{Re}}
\renewcommand{\Im}{\operatorname{Im}}

\newcommand{\bbg}[1][\Sigma]{\bg_{#1}}
\newcommand{\bbh}[1][\Sigma]{\bh_{#1}}
\newcommand{\bbt}[1][\Sigma]{\bt_{#1}}

\newcommand{\bg}{\boldsymbol{g}}
\newcommand{\bh}{\boldsymbol{h}}
\newcommand{\bt}{\boldsymbol{t}}
\newcommand{\umult}{\mathbf u}

\newcommand{\CN}{C_{\mathcal N}}
\newcommand{\tCN}{\widetilde{C}_{\mathcal N}}
\newcommand{\NormHD}[3][\Gamma]%
{\left\|#3\right\|_{H^{1/2}(#1),#2}}
\newcommand{\NormHDalt}[3][\Gamma]%
{\left\lbrace\kern-2.5pt\left\vert#3\right\vert\kern-2.5pt\right\rbrace_{H^{1/2}(#1),#2}}
\newcommand{\NormHN}[3][\Gamma]%
{\left\|#3\right\|_{H^{-1/2}(#1),#2}}
\newcommand{\NormHNalt}[3][\Gamma]%
{\left\lbrace\kern-2.5pt\left\vert#3\right\vert\kern-2.5pt\right\rbrace_{H^{-1/2}(#1),#2}}

\newcommand{\gD}{g_{\mathrm{D}}}

\newcommand{\gN}{g_{\mathrm{N}}}

\newcommand{\hD}{h_{\mathrm{D}}}
\newcommand{\hN}{h_{\mathrm{N}}}
\newcommand{\oldvarphi}{v}
\newcommand{\tracevara}{g}
\newcommand{\tracevarb}{g}
\newcommand{\slow}[1][s]{\underline{\sigma}(#1)}
\newcommand{\sLow}[1][s]{\overline{\sigma}(#1)}

\renewcommand{\d}{\;\mathrm{d}}

\newcommand{\traceD}[1][B_R]{\operatorname{\gamma}_{\mathrm{D},#1}}
\newcommand{\traceN}[1][B_R]{\operatorname{\gamma}_{\mathrm{N},#1}}
\newcommand{\traceNu}[1][\Omega]{\operatorname{\gamma}_{\nu,#1}}
\newcommand{\traceC}[1][\Sigma]{\operatorname{\boldsymbol{\gamma}}_{{#1}}}
\newcommand{\traceDext}[1][B_R]{\operatorname{\gamma}_{\mathrm{D},#1}^{\mathrm{ext}}}
\newcommand{\traceNext}[1][B_R]{\operatorname{\gamma}_{\mathrm{N},#1}^{\mathrm{ext}}}
\newcommand{\traceNuext}[1][\Omega]{\operatorname{\gamma}_{\nu,#1}^{\mathrm{ext}}}
\newcommand{\traceSD}{\traceD[B_R]}

\newcommand{\mean}[1]{\{\kern-4.2pt\{#1\}\kern-4.2pt\}}
\newcommand{\meanD}[2][B_R]{\mean{#2}_{\mathrm{D},#1}}
\newcommand{\meanN}[2][B_R]{\mean{#2}_{\mathrm{N},#1}}
\newcommand{\jump}[1]{\left[#1\right]}
\newcommand{\jumpD}[2][B_R]{\jump{#2}_{\mathrm{D},#1}}
\newcommand{\jumpN}[2][B_R]{\jump{#2}_{\mathrm{N},#1}}

\newcommand{\jumpNext}[2][B_R]{\left[#2\right]_{\mathrm{N},#1}^{\mathrm{ext},s}}

\newcommand{\jumpNextdual}[2][B_R]{\left[#2\right]_{\mathrm{N},#1}^{\mathrm{ext},\conj s}}

\theoremstyle{definition}

\theoremstyle{plain}
\newtheorem{theorem}{Theorem}[section]
\newtheorem{lemma}[theorem]{Lemma}
\newtheorem{remark}[theorem]{Remark}

\newtheorem{assumption}[theorem]{Assumption}

\newcommand{\Const}[1]{C_{\mathrm{#1}}}

\usepackage{totcount}
\newtotcounter{nconst}
\makeatletter
\newcount\rc@count
\let\rc@clearconstantlist\empty

\newcommand\rc@clearconstant[1]{\global\expandafter\let\csname rc@const@#1\endcsname\undefined}

\newcommand\resetconstants[1]{%
    \def\rc@constname{#1}
    \global\rc@count=1\relax 
    \bgroup 
        \let\\\rc@clearconstant 
        \rc@clearconstantlist
        \global\let\rc@clearconstantlist\empty 
    \egroup
}
\newcommand\const[1]{%
    \@ifundefined{rc@const@#1}{%
        \expandafter\xdef\csname rc@const@#1\endcsname{%
           \noexpand\rc@useconst{\rc@constname}{\the\rc@count}%
        }%
        \g@addto@macro\rc@clearconstantlist{\\{#1}}%
		\global\setcounter{nconst}{\the\rc@count}\relax
        \global\advance\rc@count1\relax
    }{}%
    \csname rc@const@#1\endcsname
}

\newcommand\rc@useconst[2]{\ensuremath{#1_{#2}}}
\makeatother

	\resetconstants{C}

\newcommand{\hSigma}{{\Sigma}}
\title{Stable skeleton integral equations for general coefficient Helmholtz
transmission problems
}
\author{Benedikt Gr\"a{\ss}le\footnote{
	  Institut für Mathematik, Universität Zürich, Winterthurerstr.~190, CH-8057 Zürich, Switzerland.
  (\{benedikt.graessle, stas\}@math.uzh.ch)}
  \and Ralf Hiptmair\footnote{Seminar für Angewandte Mathematik, ETH Zürich, Zürich, Switzerland.
  (ralf.hiptmair@sam.math.ethz.ch)}
	\and Stefan Sauter\footnotemark[1]
}
\date{}%
\begin{document}

\maketitle

\begin{abstract}
	A novel variational formulation of layer potentials and boundary integral operators generalizes their classical
	construction by Green's functions, which are not explicitly available for Helmholtz problems with variable coefficients.
	Wavenumber explicit estimates and properties like jump conditions follow directly from their variational definition and
	enable a non-local (``integral'') formulation of acoustic transmission problems (TP) 
	with piecewise Lipschitz coefficients.
	We obtain the well-posedness of the integral equations directly 
	from the stability of the underlying TP. The simultaneous
	analysis for general dimensions and complex wavenumbers (in this paper) 
	imposes an artificial boundary on the external Helmholtz problem and 
	employs recent insights into the associated Dirichlet-to-Neumann map.

\end{abstract}



\section{Introduction}\label{sec:Introduction}

Time-harmonic wave propagation in both homogeneous and non-homogeneous media 
is a fundamental phenomenon encountered across various scientific and engineering 
disciplines, including medical imaging, antenna design,
noise control, and radar and sonar detection.
In most practical applications, the ambient physical medium is
heterogeneous and may occupy multiple regions with distinct acoustic properties.
Typical examples are water, air, layered soil, and geological formations, each
characterized by varying propagation parameters such as density and wave speed.
These inhomogeneities introduce significant challenges in mathematical modelling,
which is crucial for improving physical understanding and enabling reliable 
numerical simulations.

The method of boundary integral equations (BIE) and their associated fast numerical 
solvers have been extensively developed for wave propagation in homogeneous media,
where they provide efficient and well-conditioned formulations.
However, the extension to heterogeneous media is non-trivial with classical 
techniques due to the presence of varying (and, in particular, non-constant) 
coefficients in the underlying partial differential equation (PDE).
This paper derives novel well-posed boundary integral formulations
for acoustic wave propagation in some media that are only required to be homogeneous
outside some bounded region and allows purely imaginary wavenumbers,
overcoming the limitations of%
~\cite{EFHS:StableBoundaryIntegral2021,FHS:SkeletonIntegralEquations2024}.
The mathematical model is the Helmholtz equation
\begin{subequations}
\label{Helm}
\end{subequations}%
\begin{equation}
-\operatorname{div}\left(  \mathbb{A}\nabla u\right)  +s^{2}pu=F\quad\text{in
}\Omega \tag{%
\ref{Helm}%
a}\label{HelmEq1}%
\end{equation}
with variable coefficients $\opA$ and $p$ on the unbounded Lipschitz domain
$\Omega\subset\mathbb{R}^{n}$.
The compact boundary $\partial\Omega$ models the surface of a scatterer and 
is partitioned into a relatively closed \emph{Dirichlet part} 
$\Gamma _{\operatorname*{D}}$ and 
a \emph{Neumann part}%
\ $\Gamma_{\operatorname*{N}}:=\partial\Omega\backslash\Gamma
_{\operatorname*{D}}$ with the boundary conditions
\begin{equation}%
\begin{aligned}
u|_{\Gamma_{\operatorname{D}}}&=g_{\operatorname*{D}} &
\text{on }\Gamma_{\operatorname*{D}},\\
(\opA\nabla u\cdot \nu)|_{\Gamma_{\operatorname{N}}}&=g_{\operatorname*{N}} &
\text{on }\Gamma_{\operatorname*{N}}
\end{aligned}
\tag{%
\ref{Helm}%
b}\label{HelmBC}%
\end{equation}
for $\gD\in H^{1/2}(\Gamma_{\mathrm{D}})$ and $\gN\in H^{-1/2}(\Gamma_{\mathrm{N}})$,
where $\nu$ denotes the outer unit normal on $\partial\Omega$.
To close the Helmholtz problem~\eqref{Helm} in the unbounded domain $\Omega$, we
impose the \emph{Sommerfeld radiation condition} towards infinity
\begin{equation}
\lim_{r\rightarrow\infty}r^{\left(  n-1\right)  /2}\left(  \partial
_{r}{u}+s{u}\right)  =0\qquad\text{with }\partial_{r}{u}=\nabla {u}\cdot\frac
{x}{\left\vert x\right\vert }\text{ uniformly in }x/\left\vert x\right\vert .
\tag{%
\ref{Helm}%
c}\label{radicond}%
\end{equation}
For plain scattering at the obstacle surface $\partial\Omega$, the source term $F$
in~\eqref{HelmEq1} is zero and the boundary data $(\gD,\gN)$ in~\eqref{HelmBC}
is given by an incident wave.
The point is that we only impose very weak conditions:
\begin{itemize}
	\myitem{\bf (C1)}\label{ass:C1} 
	The wavenumber 
	$s\in\mathbb{C}_{\geq0}^*\coloneqq\{z\in\mathbb{C}\setminus\{0\}\ :\ \Re z\geq0\}$
	has non-negative real part.
	\myitem{\bf (C2)}\label{ass:C2} 
	The coefficients 
	$\mathbb{A}\in L^{\infty}\left(\Omega;\mathbb{S}^n\right)$ and $p\in L^{\infty}\left( \Omega;\mathbb{R}\right)$
	in~\eqref{HelmEq1} 
	satisfy
	\[%
	\begin{array}[c]{cc}%
		a_{\min}\left\vert \xi\right\vert ^{2}\leq\mathbb{A}\left(  x\right)
		\mathbb{\xi}\cdot\mathbb{\xi}\leq a_{\max}\left\vert \xi\right\vert ^{2} &
		\text{for all }\xi\in\mathbb{R}^{n},\\
		p_{\min}\leq p\left(  x\right)  \leq p_{\max} &
	\end{array}
	\]
	for almost every $x\in\Omega$ with constants 
	$a_{\min}$, $p_{\min}\in\left(  0,1\right]$ and 
	$a_{\max}$, $p_{\max}\in\left[  1,\infty\right)$,
	where $\mathbb{S}^n$ denotes the symmetric $n\times n$ matrices.
	\myitem{\bf (C3)}\label{ass:C3}
	There is an open ball $B_{R}$ of sufficiently large radius $R>0$ 
	about the origin with
	\[
	\operatorname*{supp}\left(  \mathbb{I}-\mathbb{A}\right)  \cup
	\operatorname*{supp}\left(  1-p\right)  \cup\left(  \mathbb{R}^{n}%
	\backslash\Omega\right)  \subset
	B_{R},
	\]
	where $\mathbb{I}$ is the identity matrix and $1$ denotes the constant
	function one.
	\myitem{\bf (C4)}\label{ass:C4}
	The volume source $F$ is supported in the closure of the ball $B_R$ from~\ref{ass:C3},
	i.e.,
	\begin{align*}
		\mathrm{supp}(F)\subset\overline{B_R}.
	\end{align*}

\end{itemize}
The conditions \ref{ass:C3}--\ref{ass:C4} are classical for exterior Helmholtz problems with variable
coefficients~\cite{BCT:UniquenessTimeHarmonic2012,GPS:HelmholtzEquationHeterogeneous2019,SW:WavenumberexplicitParametricHolomorphy2023}
and imply that the Helmholtz equation~\eqref{HelmEq1} 
has a homogeneous far-field.

Conditions~\ref{ass:C2}--\ref{ass:C3} permit a broad class of coefficients 
in the Helmholtz equation, specifically allowing for applications with piecewise smooth
material parameters $\mathbb{A}$ and $p$, which represent,
e.g., different media with varying properties inside a large ball.
\begin{figure}[t]
	\centering
	\includegraphics[]{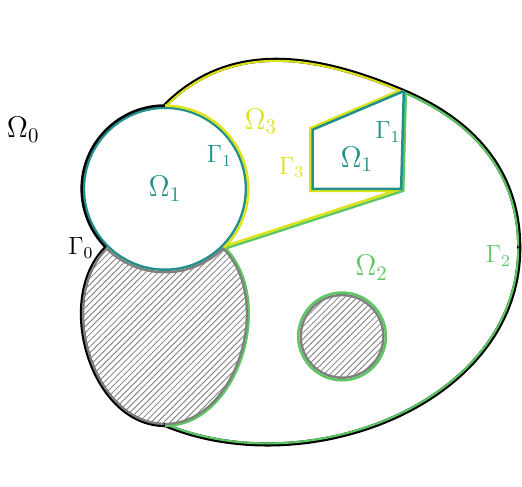}
	\caption{Decomposition of $\Omega$ into Lipschitz sets 
		$\Omega_0,\dots,\Omega_J$ with their respective boundaries 
		$\Gamma_0,\dots,\Gamma_J$ for $J=3$ and the acoustic obstacle $\R^\dim\setminus\Omega$ (hatched grey) in
\Cref{sec:Stable integral formulation for the transmission problem}.}
	\label{fig:acoustic_domain}
\end{figure}%
This suggests a decomposition of $\Omega$ into $J\in\mathbb{N}$
disjoint, open, and bounded Lipschitz sets $\Omega_{1},\dots,\Omega_J\subset\Omega$
and the unbounded complement
\[
\Omega_{0}\coloneqq\Omega\backslash\left(  \bigcup_{j=1}^{J}\overline
{\Omega_{j}}\right)
\]
as illustrated in \Cref{fig:acoustic_domain},
such that the restrictions $\opA_j\coloneqq\left.  \mathbb{A}%
\right\vert _{\Omega_{j}}$ and $p_j\coloneqq\left.  p\right\vert _{\Omega_{j}}$ 
are smooth or even constant. 
In the case of pure scattering problems ($F=0$), the original problem~\eqref{Helm}
can be formulated as a Helmholtz transmission problem for the solution $u_j\coloneqq u|_{\Omega_j}$ 
over the decomposition into $\Omega_j$ with outer unit normal $\nu_j$, namely
\begin{subequations}\label{eqn:ATP_intro}
	\begin{align}
		-\Div(\opA_j \nabla u_j) +s^2 p_j\, u
			&= 0 
			&&\text{in }\Omega_j\text{ for }j=0,\dots,J,\label{eqn:ATP_intro_a}\\
		(\opA_j\nabla u_j\cdot \nu_j)|_{\Gamma_j\cap \Gamma_k} 
				+(\opA_k\nabla u_k\cdot\nu_k)|_{\Gamma_j\cap \Gamma_k}
				&= 0
				&&\text{on }\Gamma_j\cap\Gamma_k\text{ for }j,k=0,\dots,J,\label{eqn:ATP_intro_b}\\
		u_j|_{\Gamma_j\cap \Gamma_k}-u_k|_{\Gamma_j\cap \Gamma_k} &= 0
				&&\text{on }\Gamma_j\cap\Gamma_k\text{ for }j,k=0,\dots,J,\label{eqn:ATP_intro_c}\\
				(\opA_j\nabla u_j\cdot \nu_j)|_{\Gamma_j\cap \Gamma_{\mathrm{N}}} 
				&= \gN|_{\Gamma_{\mathrm{N}}}
				&&\text{on }\Gamma_j\cap\Gamma_{\mathrm{N}}\text{ for }j=0,\dots,J,\label{eqn:ATP_intro_d}\\
		u_j|_{\Gamma_j\cap \Gamma_{\mathrm{D}}}&=\gD|_{\Gamma_{\mathrm{D}}}	 
			   &&\text{on }\Gamma_j\cap\Gamma_{\mathrm{D}}\text{ for }j=0,\dots,J,\label{eqn:ATP_intro_e}\\
			   \span\span u_0\text{ satisfies }~\eqref{radicond}.\span\span\label{eqn:ATP_intro_f}
	\end{align}
\end{subequations}
This paper presents the method of
\textit{skeleton integral equations} (SIE) to transform~\eqref{eqn:ATP_intro}
for general coefficients and wavenumbers
in a natural way to a non-local (integral%
\footnote{Throughout this paper, we employ the term ``integral'' equations 
	for non-local operator equations. In the case of certain ``nice'' (e.g., constant) 
	coefficients $\opA$ and $p$, these operators have classical 
	integral representations with known kernel functions.
}%
) equation on the 
\textit{skeleton} $\bigcup_{j=0}^{J}\partial\Omega_{j}$ such that our main paradigm applies:
\begin{align}\label{eqn:main_par}
	\parbox{0.8\linewidth}{
	\textquotedblleft The Helmholtz problem~\eqref{Helm} is well posed
	if and only if the skeleton integral equation is well posed.\textquotedblright }
\end{align}

\textbf{Main results and literature review}. 
The main task in deriving the SIE in the present setting is
the generalisation of the classical layer potentials, which are based on the explicit knowledge of 
the Green's function that is not available here.
Earlier works~\cite{FHS:SkeletonIntegralEquations2024,EFHS:StableBoundaryIntegral2021}
on general coefficients $\opA$ and $p$ with \ref{ass:C2}
introduced a variational definition of layer potentials 
for wavenumbers with positive real part.
The analysis therein relies on the $H^1(\Omega)$-coercivity of the sesquilinear form
for the variational Helmholtz problem~\eqref{Helm} and
breaks down as the real part of the wavenumber tends to zero. 
This paper presents a unified analysis for wavenumbers $s\in\mathbb{C}_{\geq0}^*$ 
and generalises the
\textit{multi-trace} and \textit{single-trace formulations} introduced in
\cite{CH:IntegralEquationsAcoustic2015,CHJP:NovelMultitraceBoundary2015}
to Helmholtz transmission
problems with varying coefficients (beyond the case of piecewise constants).
For an overview of
the many ways to transform the PDEs \eqref{Helm} to integral equations
on the domain skeleton, we refer to%
~\cite{BLS:StableNumericalCoupling2015,Say:RetardedPotentialsTime2016} and the references therein.
The main methodological and theoretical results are summarised as follows.

\begin{enumerate}[label=(\alph*)]
\item 
In the case of purely imaginary wavenumbers $s\in \operatorname{i}\R$ in $n=2,3$ dimensions and constant (isotropic) media,
it is a classical approach%
~\cite{Ned:AcousticElectromagneticEquations2001,MS:ConvergenceAnalysisFinite2010,GPS:HelmholtzEquationHeterogeneous2019}
to consider an equivalent reformulation of the indefinite Helmholtz equation~\eqref{HelmEq1}
on the finite domain $B_R\cap \Omega$ with
\emph{Dirichlet-to-Neumann}
($\operatorname*{DtN}$) \emph{boundary conditions} on the artificial boundary $\partial B_R$.
Our unified analysis extends this technique
to general wavenumbers $s\in\mathbb{C}_{\geq0}$ and spatial dimensions $n\geq2$ 
based on new results for the DtN operator in%
~\cite{GS:DirichlettoNeumannOperatorHelmholtz2025}.

\item For problems with constant coefficients, it is well known that
Helmholtz-harmonic functions on bounded domains can be represented by means of
their Cauchy traces on the boundary%
~\cite[Thm.~3.1.8]{SS:BoundaryElementMethods2011} 
based on the single and double layer
potential operators in a Green's representation formula. For our setting with
$L^{\infty}$ coefficients $\mathbb{A}$ and $p$ and (possibly) purely imaginary wavenumbers $s\in
\operatorname*{i}\mathbb{R}$, the standard definitions in the literature are
not applicable. We define these layer operators in this paper for the
class of coefficients satisfying \ref{ass:C2}--\ref{ass:C3} and 
establish their defining mapping properties and jump conditions.
From this, we obtain Green's representation formula and deduce the
Calder\'{o}n identity for the Cauchy traces. 
The layer potentials will be defined as solutions of certain transmission problems
and we prove their well-posedness even for the critical case of purely
imaginary wavenumbers.

\item The Cauchy data of the transmission problem~\eqref{eqn:ATP_intro}
	satisfy Calder\'{o}n identities for each subdomain
boundary and are subject to the transmission and boundary conditions. This leads to a
\emph{multi-trace formulation} of the Helmholtz transmission problem. 
An equivalent formulation for classical \emph{single-trace spaces} with incorporated interface
and boundary conditions results in the \emph{single-trace formulation} of~\eqref{eqn:ATP_intro}.
Our main paradigm~\eqref{eqn:main_par} implies
the well-posedness of these novel \emph{skeleton integral equations}
from that of the original Helmholtz problem~\eqref{Helm}.
\end{enumerate}

This paper is a contribution to the analysis of general Helmholtz PDE and
provides wavenumber-explicit stability estimates.
We establish that the well-posedness of the
novel SIE introduced in this paper is unconditionally equivalent
to the well-posedness of the original PDE. This generalises the
method of integral equations for PDE from constant coefficients
to variable, rough coefficients and provides the theoretical
tools for their analysis. These SIE serve as a starting
point for numerical discretisations, e.g., by the boundary
element method.
In particular, the modelling of \textit{high-frequency} scattering problems in heterogeneous media
by integral equations with non-local $\operatorname*{DtN}$ boundary conditions
is appealing for various practical reasons; among the most important are:

\begin{enumerate}[label=(\Alph*)]
\item Highly indefinite Helmholtz PDEs with \textit{approximate}
local or non-local boundary conditions of, e.g., impedance/Robin type
\cite{Gol:FiniteElementMethod1982}, non-local absorbing boundary conditions%
~\cite{Tsy:NumericalSolutionProblems1998,Giv:NumericalMethodsProblems1992,HMG:HighorderLocalAbsorbing2008,Ihl:FiniteElementAnalysis1998}
or perfectly matched
layers%
~\cite{Ber:PerfectlyMatchedLayer1994,ST:PerfectlyMatchedLayer2004,BBL:PerfectlyMatchedLayers2003}
may lead to significant pollution and stability
issues.
In contrast, the $\operatorname*{DtN}$ condition allows exact representations 
(without further approximation) by boundary integral operators%
~\cite[p.~52]{CK:InverseAcousticElectromagnetic2019}

\item The method of integral equations reduces the Helmholtz PDE in $n$
spatial dimensions to the $\left(  n-1\right)$-dimensional domain skeleton,
whose numerical discretisation requires fewer degrees of freedom for similar accuracies
(see, e.g.,~\cite{SS:BoundaryElementMethods2011}). 
Our key paradigm~\eqref{eqn:main_par} reduces the well-posedness to standard results in the literature.
The well-posedness of \eqref{Helm} is
established for piecewise constant coefficients%
~\cite{CK:InverseAcousticElectromagnetic2019,KR:TransmissionProblemsHelmholtz1978,McL:StronglyEllipticSystems2000,Von:BoundaryIntegralEquations1989}
and follows from~\cite{BCT:UniquenessTimeHarmonic2012} 
for piecewise Lipschitz coefficients, see, e.g.,%
~\cite{GPS:HelmholtzEquationHeterogeneous2019,SW:WavenumberexplicitParametricHolomorphy2023}.
Frequency-explicit estimates exist if the matrix coefficient 
$\mathbb{A}$ in \eqref{HelmEq1} satisfies
certain monotonicity or regularity assumptions%
~\cite{Bur:DecroissanceLenergieLocale1998,EM:StabilityDiscretizationsHelmholtz2012,HPV:AsymptoticallyPreciseNorm2007,BCT:UniquenessTimeHarmonic2012,MS:HelmholtzEquationReally2014,GPS:HelmholtzEquationHeterogeneous2019,GS:StabilityFiniteElement2019,ST:HeterogeneousHelmholtzProblem2021,GSW:OptimalConstantsNontrapping2020}.
\end{enumerate}

\noindent\textbf{Outline and further contributions.}
\Cref{sub:general_notation} introduces the geometric setting and the general notation of
Sobolev spaces and their traces.

The Helmholtz equation~\eqref{Helm} on the unbounded domain $\Omega$ and its 
equivalent strong and variational formulations on the truncated domain 
$\Omega\cap B_R$ are discussed with the appropriate Sobolev setting for the given data $\gD,\gN, F$ and solution $u$ in \Cref{sec:potential_operators}.
A Fredholm argument provides the equivalence of the well-posedness of the Helmholtz problem
and the uniqueness of its solutions in \Cref{thm:L_N_well_def}.
The uniqueness may follow from a unique continuation principle for piecewise Lipschitz
coefficient matrix $\opA$~\cite{BCT:UniquenessTimeHarmonic2012} 
and is supposed throughout 
this paper by \Cref{ass:A_p}.
\Cref{sub:Newton_potential} investigates the continuous solution operator
for the truncated domain.
Particular care is taken in \Cref{lem:Newton_potential_L2} 
to characterise its restriction to more regular $L^2$ sources
which significantly simplifies part of the following analysis.

Section \ref{sec:Potential operators for interface problems} 
defines layer potentials for the class of general coefficients with \ref{ass:C2}--\ref{ass:C3} 
as solutions to variational transmission problems in the full space $\mathbb{R}^{n}$.
The definition of the single layer operator and the verification of its
jump conditions extends
\cite[Lem.~3.7]{FHS:SkeletonIntegralEquations2024}.
We establish a natural operator representation of the single and double layer potentials
as the composition of the solution operator with dual trace operators
that generalise their definitions%
~\cite{McL:StronglyEllipticSystems2000,SS:BoundaryElementMethods2011} 
in the homogeneous case and was only known
for the single layer operator%
~\cite{Bar:LayerPotentialsGeneral2017,FHS:SkeletonIntegralEquations2024}.
For the double layer potential, the representation
relies on an extension of the Newton potential that was not available before.
This obstacle motivated alternative definitions of the double layer potential in%
~\cite[Sec.~4-5]{Bar:LayerPotentialsGeneral2017}
based on a lifting of the boundary density 
and, for the definite case, in~\cite[Subsec.~3.2.3]{FHS:SkeletonIntegralEquations2024} 
through a variational formulation
that is generalised by our natural definition.
The analysis of the double layer potentials is more involved and 
requires a detour over a mixed reformulation 
to prove its mapping properties and jump conditions.
This enables a Green's representation
formula and the application of the Cauchy trace operator leads to the
Calder\'{o}n identities in \Cref{sub:The Calder'on operator}.

With these definitions at hand, we derive in \Cref{sec:Stable integral formulation for the transmission problem} stable
SIE formulation of our transmission problem; first, in a multi-trace setting with
the transmission and boundary condition as additional constraints and then as
a single-trace integral equation in operator form. 
Theorem \ref{thm:TP_equivalence} proves our main paradigm~\eqref{eqn:main_par}
on the equivalence of the original Helmholtz problem (\ref{Helm}) and
the multi- and single-trace formulations.

\section{Preliminaries}%
\label{sub:general_notation}
A domain is a (possibly unbounded) nonempty, open, and connected subset $\omega\subset\R^n$ of the $n$-dimensional Euclidean space. 
It is an exterior domain if its complement $\R^n\setminus\omega$ is bounded.
The set of non-zero complex numbers 
with non-negative real part reads
\begin{align*}
	\Cstar_{\geq0}\coloneqq\{z\in \mathbb{C}\ :\ \Re z\geq0
	\text{ and }z\ne0\}.
\end{align*}
Standard notation on (complex-valued) Lebesgue and Sobolev spaces and their norms
applies for open subsets $\omega\subset\R^n$ with $n\geq2$ 
throughout this paper.
In particular, the space $H_{\operatorname*{loc}%
}^{\kappa}\left(  \omega\right)  $ for $\kappa\geq0$
is given by all distributions $v\in (C_{\mathrm{comp}}^\infty(\omega))'$
on the compactly supported smooth functions 
$C_{\operatorname*{comp}}^{\infty}\left(  \omega\right)  $ such that $\varphi v\in H^{\kappa}\left(
\omega\right)  $ for all $\varphi\in C_{\operatorname*{comp}}^{\infty}\left(
\R^n\right)  $.
For the definition of Sobolev spaces $H^{\kappa}\left(
\Gamma\right)  $ on relatively open parts $\Gamma\subset \partial\omega$
of Lipschitz boundaries $\partial\omega$
and their norm $\left\Vert \cdot\right\Vert _{H^{\kappa}\left(
\Gamma\right)  }$ we refer, e.g., to~\cite[pp.~96--99]{McL:StronglyEllipticSystems2000}.

The natural energy norms for $H^1(\omega)$ and $H(\omega,\Div)$ 
for the wavenumber $s\in\Cstar_{\geq0}$ read
\begin{align}\label{eqn:H_norm_def}
	\|v\|_{H^1(\omega),s}
	&\coloneqq\sqrt{\|\nabla v\|_{L^2(\omega)}^2 + |s|^2\|v\|_{L^2(\omega)}^2}
	&&\text{for all }v\in H^1(\omega),\\
	\|\mathbf{p}\|_{H(\omega,\Div),s}
	&\coloneqq\sqrt{|s|^{-2}\|\Div \mathbf{p}\|_{L^2(\omega)}^2 
	+ \|\mathbf{p}\|_{L^2(\omega)}^2}
	&&\text{for all }\mathbf{p}\in H(\omega,\Div)\label{eqn:Hdiv_norm_def}.
\end{align}
The anti-dual version of the $L^2$ scalar product on $L^2(G)$
for an open Lipschitz set $G=\omega\subset\R^\dim$ with outer unit normal $\nu_{\omega}$
or its boundary $G=\partial\omega$ is written as
\begin{align*}
	\left\langle v,w\right\rangle_{G}
	\coloneqq\int_{G} v\, w\d x(G) 
	\qquad\text{for all }v,w\in L^2(G),
\end{align*}
and extends to the natural dual pairing on 
$H^{1/2}(G)\times H^{-1/2}(G)$
(and $H^{-1/2}(G)\times H^{1/2}(G)$) with the same notation.
The (Dirichlet) trace operator $\traceD[\omega]:H^1(\omega)\to H^{1/2}(\partial\omega)$ 
is surjective and the unique continuous operator with
\begin{align*}
	\traceD[\omega] v = v|_{\partial\omega}\qquad\text{for all }
	v\in C^\infty(\overline{\omega}).
\end{align*}
The topological dual space $H^{-1/2}(\partial\omega)=(H^{1/2}(\partial\omega))'$ 
consists of all normal traces 
$\traceNu[\omega]\mathbf{q}\coloneqq \mathbf{q}|_{\partial\omega}\cdot \nu_{\omega}$ 
of functions $\mathbf{q}\in H(\omega,\Div)$ that are defined by
\begin{align*}
	\left\langle \traceNu[\omega]\mathbf{q}, \traceD[\omega] v\right\rangle_{
			\partial\omega} = \int_\omega
	v \Div \mathbf{q} + \mathbf{q}\cdot \nabla v\d x\qquad\text{for all }v\in H^1(\omega).
\end{align*}
The associated Sobolev spaces with boundary conditions on $\Gamma\subset\partial\omega$ read
\begin{align*}
	H^1_{\Gamma}(\omega)&\coloneqq\{v\in H^1(\omega)\ :\ v|_{\Gamma}\equiv 0 \},\\
	H_{\Gamma}(\omega,\Div)&\coloneqq\{\mathbf{q}\in H(\omega,\Div)\ :\ (\mathbf{q}\cdot\nu_{\omega})|_{\Gamma}\equiv 0
	\}.
\end{align*}
Given $\opA\in L^\infty(\omega;\mathbb{S}^n)$ 
with values in the symmetric $n\times n$ matrices
$\mathbb{S}^n\subset\R^{n\times n}$,
the spaces
\begin{align*}
	H^1(\omega,\opA)&\coloneqq\{ v\in H^1(\omega)\ :\ \opA\nabla v\in
	H(\omega,\Div)\},\\
	H^1_{\mathrm{loc}}(\omega,\opA)&\coloneqq\{ v\in H^1_{\mathrm{loc}}(\omega)\ 
		:\ \varphi\opA\nabla v\in
	H(\omega,\Div)\text{ for all }\varphi\in C^\infty_{\mathrm{comp}}(\R^n)\}
\end{align*}
admit a continuous trace operator 
$\traceN[\omega]:H^1_{\mathrm{loc}}(\omega,\opA)\to H^{-1/2}(\partial\omega)$ with
\begin{align*}
	\traceN[\omega] v\coloneqq \traceNu[\omega](\opA\nabla v)
	\qquad\text{for all }
	v\in H^1_{\mathrm{loc}}(\omega,\opA),
\end{align*}
called the (co-)normal (or Neumann) trace operator.
The dependence of $\traceN[\omega]{}$ on the matrix function $\opA$ 
will be clear from the context and is surpressed in this notation.
The identity matrix is denoted by $\mathbb{I}\subset\R^{\dim\times \dim}$.
The trace operators $\traceD[\omega], \traceN[\omega]$, 
and $\traceNu[\omega]$ may also be applied to
appropriate Sobolev functions defined on another Lipschitz set
$\omega_0\subset\R^\dim$, whose closure contains 
$\partial\omega\subset\overline{\omega_0}$, 
and always denotes the corresponding trace on $\partial\omega$.
The Dirichlet trace $\traceDext[\omega]=\traceD[\omega_{\mathrm{ext}}]$ 
and Neumann trace $\traceNext[\omega]=\traceN[\omega_{\mathrm{ext}}]$ 
on the exterior domain 
$\omega_{\mathrm{ext}}\coloneqq\R^\dim\setminus\overline{\omega}$ with outer normal 
$\nu_{\omega_{\mathrm{ext}}}=-\nu_\omega$ on $\partial\omega$ define the jumps
and averages by
\begin{equation}\label{eqn:jump_mean_def}
	\begin{aligned}
		\jumpD[\omega]{v}
		&\coloneqq \traceD[\omega] v - \traceDext[\omega] v,
		\quad
		\meanD[\omega]{v}\coloneqq \tfrac{1}{2}(\traceD[\omega]v+\traceDext[\omega]v)
		&&\text{for all }v\in
		H^1(\R^\dim\setminus\partial\omega),\\
		\jumpN[\omega]{v}
		&\coloneqq\traceN[\omega] v + \traceNext[\omega] v,
		\quad
		\meanN[\omega]{v}\coloneqq \tfrac{1}{2}(\traceN[\omega]v-\traceNext[\omega]v)
		&&\text{for all }v\in H^1(\R^\dim\setminus\partial\omega,\opA).
	\end{aligned}
	\end{equation}
The open ball of radius $R>0$ about the origin is denoted by
\begin{align*}
	B_R\coloneqq\{x\in \R^n\ :\ \|x\|<R\}
\end{align*}
with Euclidean norm $\|\bullet\|$.
Its outer unit normal vector $x/\left\Vert x\right\Vert$ 
on the boundary $S_{R}=\partial B_R$ points into the unbounded complement 
$B_R^+:=\mathbb{R}^n\setminus\overline{B_R}$ of $B_{R}$ and the normal derivative in this direction 
is denoted by $\partial_{r}$. 
The notation $|\bullet|$ is context-sensitive and may refer to the Lebesgue measure $|\omega|$ of a
bounded measurable $n$-dimensional set $\omega\subset\mathbb{R}^{n}$, 
the surface measure $|\Gamma|$ of an $\left(
n-1\right)  $-dimensional manifold $\Gamma\subset\mathbb{R}^{n}$, 
the cardinality $|J|$ of a countable set $J$, and the absolute value $|z|$
of a complex number $z\in\mathbb{C}$.

\section{Helmholtz problem with varying coefficients}%
\label{sec:potential_operators}
The well-posedness of the exterior Dirichlet
problem outside a large ball 
enables a unified analysis of the exterior Helmholtz problem for general coefficients and
wavenumbers $s\in\Cstar_{\geq0}$.

\subsection{The exterior Helmholtz problem}%
\label{sub:The full space Helmholtz problem}
Let $\Omega\subset\R^\dim$ denote an unbounded Lipschitz domain 
in $\dim\geq2$ dimensions with bounded 
(and possibly multiply connected or empty) complement $\R^\dim\setminus\Omega$.
The compact boundary $\partial\Omega$
splits into a relatively closed Dirichlet part $\Gamma_{\mathrm{D}}$ and 
the Neumann part $\Gamma_{\mathrm{N}}=\partial\Omega\setminus\Gamma_{\mathrm{D}}$.
The space
\begin{align*}
	H^{-1}_{\mathrm{D},\mathrm{comp}}(\Omega)
	\coloneqq\left\{F\in (H^1_{\Gamma_\mathrm{D}}(\Omega))' \ :\ \mathrm{supp}(F)\text{ is compact}\right\},
\end{align*}
contains the admissible sources with compact support (relative to $\R^\dim$)
in the dual space of 
$H^1_{\Gamma_{\mathrm{D}}}(\Omega)$.
The exterior Helmholtz problem~\eqref{Helm}
with source $F\in H^{-1}_{\mathrm{D},\mathrm{comp}}(\Omega)$,
Dirichlet data $\gD\in H^{1/2}(\Gamma_{\mathrm{D}})$,
and Neumann data $\gN\in H^{-1/2}(\Gamma_{\mathrm{N}})$
seeks a solution $u\in H^1_{\mathrm{loc}}(\Omega)$ to
\begin{equation}\label{eqn:FSP}
	\begin{aligned}
		-\Div(\opA\nabla u) + s^2 p u 
		&= F &&\text{in }\Omega,\\
		u|_{\Gamma_{\operatorname{D}}}&=g_{\operatorname*{D}} &&
\text{on }\Gamma_{\operatorname*{D}},\\
		(\opA\nabla u\cdot \nu)|_{\Gamma_{\operatorname{N}}}&=g_{\operatorname*{N}} &&
\text{on }\Gamma_{\operatorname*{N}},\\
		u&\text{ satisfies }\eqref{radicond}\span\span\span
	\end{aligned}
\end{equation}
for a wavenumber $s\in\Cstar_{\geq0}$, (possibly non-constant)
coefficient functions $\opA\in L^\infty(\Omega; \mathbb{S}^{\dim})$ 
and $p\in L^\infty(\Omega)$, and $F$
satisfying~\ref{ass:C1}--\ref{ass:C4} for a sufficiently large ball%
\footnote{In the coercive case $\Re s>0$, the further analysis also applies to
$R=\infty$
as discussed in \Cref{rem:B_infty}.}
$B_R$.
The point of~\ref{ass:C3}--\ref{ass:C4} is the existence 
\cite[Chap.~9]{McL:StronglyEllipticSystems2000} of a unique 
solution $u_{\mathrm{ext}}\in H^1_{\mathrm{loc}}(B_R^+,\IMatrix)$ 
to the corresponding exterior Helmholtz problem 
in $B_R^+\coloneqq \R^\dim\setminus \overline{B_R}$
for any Dirichlet data $\hD\in H^{1/2}(S_R)$
or any Neumann data $\hN\in H^{-1/2}(S_R)$ on the sphere $S_R\coloneqq\partial B_R$, 
namely
\begin{equation}\label{eqn:farfield_problem}
	\begin{aligned}
		-\Delta u_{\mathrm{ext}} + s^2 u_{\mathrm{ext}} 
			&= 0	
			&&\text{in }B_R^+,\\
		\text{either }u_{\mathrm{ext}} = \hD 
		\text{ or }
		\partial_r u_{\mathrm{ext}} &= \hN
									&&\text{on }S_R,\\
		u_{\mathrm{ext}}&\text{ satisfies }\eqref{radicond}.\span\span\span
	\end{aligned}
\end{equation}
The well-posedness~\cite{McL:StronglyEllipticSystems2000}
of~\eqref{eqn:farfield_problem} induces either of the solution maps
\begin{align*}
	\opSD&:H^{1/2}(S_R)\to H^1_{\mathrm{loc}}(B_R^+,\IMatrix)
		 &&\text{with }\opSD\hD\coloneqq u_{\mathrm{ext}}\quad\text{or }\\
	\opSN&:H^{-1/2}(S_R)\to H^1_{\mathrm{loc}}(B_R^+,\IMatrix)
		 &&\text{with }\opSN\hN\coloneqq u_{\mathrm{ext}}.
\end{align*}
Clearly, $\traceDext\circ \opSD=\id$ and $\traceNext\circ \opSN=-\id$,
where $-\traceNext = (\partial_r\bullet)|_{S_R}$ corresponds to the 
normal derivative with respect to the outer unit normal for $B_R$.
Their other traces 
define the 
\emph{Dirichlet-to-Neumann} operator
\begin{align}\label{eqn:DtN_def}
	\DtN\coloneqq -\traceNext\circ \opSD:H^{1/2}(S_R)\to H^{-1/2}(S_R)
\end{align}
and the \emph{Neumann-to-Dirichlet} operator 
\begin{align*}
	\NtD\coloneqq +\traceDext\circ \opSN:H^{-1/2}(S_R)\to H^{1/2}(S_R).
\end{align*}
The $\mathrm{DtN}$ and $\mathrm{NtD}$ operator are naturally inverse to each other.
If no confusion arises, we abbreivate here and in the remaining parts of this paper
\begin{align*}
	\opSD v\coloneqq\opSD\traceD v,
	\ \DtN v\coloneqq\DtN\traceD v,
	\ %
	\NtD v\coloneqq\NtD\traceN v.
\end{align*}
These operators enable an equivalent reformulation of
the Helmholtz problem~\eqref{eqn:FSP}
on the truncated domain $\Omega_R\coloneqq B_R\cap \Omega$ that goes back at least to
\cite{MM:FiniteElementMethod1980,Mas:NumericalSolutionExterior1987,KM:ConvergenceAnalysisCoupled1990,Ned:AcousticElectromagneticEquations2001}
for $\Re s=0$.
The resulting truncated Helmholtz problem seeks a solution $u_R\in H^1(\Omega_R)$~to 
\begin{equation}\label{eqn:THP}
	\begin{aligned}
		-\Div(\opA\nabla u_R) + s^2 p \, u_R &= F &&\text{in }\Omega_R,\\
		\partial_r u_R&=\DtN u_R &&\text{on }S_R,\\
		(\traceD[\Omega] u_R)|_{\Gamma_{\mathrm{D}}}
					  &=\gD&&\text{on }\Gamma_{\mathrm{D}},\\ 
		(\traceN[\Omega] u_R)|_{\Gamma_{\mathrm{N}}}
					  &=\gN&&\text{on }\Gamma_{\mathrm{N}}.
	\end{aligned}
\end{equation}
(The boundary condition on $S_R$ in~\eqref{eqn:THP} may be equivalently replaced by $u_R=\NtD u_R$.)
The problem~\eqref{eqn:THP} trades an additional boundary condition at an %
artificial boundary $S_R$ for the boundedness of $\Omega_R$.
Series representations and properties of $\DtN$ known from
\cite{Ned:AcousticElectromagneticEquations2001,MS:ConvergenceAnalysisFinite2010}
for $\Re s=0$ and $\dim=2,3$ are discussed in~\cite{GS:DirichlettoNeumannOperatorHelmholtz2025} 
for $s\in\Cstar_{\geq0}$ and~$\dim\geq2$.%
\begin{theorem}[{equivalence}]
	\label{lem:equivalence}
	If \ref{ass:C1}--\ref{ass:C4} hold,
	$u\in H^1_{\mathrm{loc}}(\Omega)$ is a solution to%
	~\eqref{eqn:FSP} if and only if %
	$u_R\coloneqq u|_{\Omega_R}\in H^1(\Omega_R)$ 
	solves%
	~\eqref{eqn:THP} and $u|_{B_R^+}=S_{D} u_R$.
\end{theorem}
\begin{proof}
	Any solution $u\in H^1_{\mathrm{loc}}(\Omega)$ 
	to~\eqref{eqn:FSP} with $F|_{B_R^+}\equiv 0$ satisfies 
	$u|_{B_R^+}=\opSD u$ by the uniqueness of solutions
	\cite[Thm.~9.11]{McL:StronglyEllipticSystems2000}
	of~\eqref{eqn:farfield_problem}. 
	Hence $u|_{B_R}$ satisfies~\eqref{eqn:THP}. 
	Conversely, any solution $u_R\in H^1(\Omega_R)$ 
	to~\eqref{eqn:THP} extends to $u\in H^1_{\mathrm{loc}}(\Omega)$ by
	$u|_{B_R^+}\coloneqq\opSD u_R$. 
	This and and the continuity $\partial_ru=\DtN u$ at $S_R$
	reveal~\eqref{eqn:FSP}.
\end{proof}
\subsection{Uniqueness and existence of solutions}%
\label{sub:Variational_form_existence}
In the case of absorption ($\Re s>0$),
the uniqueness and existence of solutions to~\eqref{eqn:FSP} 
-- and equivalently to~\eqref{eqn:THP} by \Cref{lem:equivalence} --
is a consequence of the continuity and coercivity 
of the associated bilinear form%
~\cite[Lem.~3.2]{FHS:SkeletonIntegralEquations2024}.
This is different for the indefinite case with $\Re s=0$,
where the well-posedness of the (truncated) Helmholtz problem 
classically follows from a Fredholm alternative argument and the
uniqueness of solutions.
\begin{assumption}\label{ass:A_p}
	For any $F\in H^{-1}_{\mathrm{D, comp}}(\Omega)$ with~\ref{ass:C4}, there exists at most one
	solution to the truncated problem~\eqref{eqn:THP}.
\end{assumption}

\Cref{ass:A_p} can be understood as an additional condition on the coefficient 
$\opA$ in the case $\Re s=0$ and holds for a large class of 
\emph{well-behaved} coefficients:
The seminal paper%
~\cite{BCT:UniquenessTimeHarmonic2012}
and~\cite[Prop.~2.13]{LRX:MoscoConvergenceHcurl2019} establish a unique continuation
principle for \emph{piecewise} Lipschitz $\opA$.
Even though those references consider Maxwell's equation, their arguments
apply to the Helmholtz equation in any dimension $\dim\geq 2$ based on the 
unique continuation property for globally Lipschitz coefficients $\opA\in
W^{1,\infty}(\R^\dim)$~\cite{AKS:UniqueContinuationTheorem1962,Wol:PropertyMeasuresRN1992}.%
\begin{lemma}[{uniqueness for piecewise Lipschitz $\opA$}]\label{lem:uniqueness}
	If there is a finite collection $(\omega_j)_{j=1}^N$ of $N\in\mathbb N$ 
	pairwise disjoint domains $\omega_j\subset\R^\dim$ 
	of class $C^0$ with $\R^\dim=\cup_{j=1}^N \overline{\omega_j}$
	and $\opA|_{\omega_j}=\opA_j$ for some
	$\opA_j\in W^{1,\infty}(\R^\dim;\mathbb S^\dim)$ and all $j=1,\dots,N$,
	then \Cref{ass:A_p} holds.
\end{lemma}
The proof of~\Cref{lem:uniqueness} utilises the sign properties of the 
DtN operator
\begin{align}\label{eqn:DtN_non_pos}
	0&\leq -\Re\left(\left\langle \DtN \tracevarb,\conj{\tracevarb}\right\rangle_{S_R}\right)
	 &&\text{for all }\tracevarb\in H^{1/2}(S_R),\\\label{eqn:DtN_Im_sign}
	0&< -\Im(s)\Im\left(\left\langle \DtN \tracevarb,\conj{\tracevarb}\right\rangle_{S_R}\right)
	 &&\text{for all }\tracevarb\in H^{1/2}(S_R)\setminus\{0\}\text{ and }\Im s\ne0
\end{align}
known from \cite{Ned:AcousticElectromagneticEquations2001,MS:ConvergenceAnalysisFinite2010} for $\Re s=0$, $\dim=2,3$
and from~\cite{GS:DirichlettoNeumannOperatorHelmholtz2025} in the general case.
\begin{proof}[Proof of \Cref{lem:uniqueness}]
	It suffices to prove that the homogeneous problem has at most one solution.
	Let $u\in H^1(\Omega_R)$ solve~\eqref{eqn:THP} 
	with vanishing data $F,\gD$, and $\gN$.
	A standard argument with an integration by parts provides
	\begin{align*}%
		\|\opA^{1/2}\nabla u\|_{L^2(\Omega_R)}^2
		+s^2\|p^{1/2}u\|_{L^2(\Omega_R)}^2
		-\left\langle \DtN u,\conj u\right\rangle_{S_R}=0.
	\end{align*}
	The multiplication in $\mathbb{C}$, the sign~\eqref{eqn:DtN_non_pos} of the real part of $\DtN$, 
	and $\Re(s)\geq0$ reveal
	\begin{align*}
		\Re\big(\conj s\left\langle \DtN u,\conj u\right\rangle_{S_R}\big)
		\leq\Im(s)\Im\big(\left\langle \DtN u, \conj{u}\right\rangle_{S_R}\big).
	\end{align*}
	This and 
	the real part of the previous identity %
	multiplied by $\conj s$ verify
	\begin{align}\label{eqn:re_norm_u_esimate}
		\Re(s)\|\opA^{1/2}\nabla u\|_{L^2(\Omega_R)}^2 
		+ \Re(s)|s|^2\|p^{1/2}u\|_{L^2(\Omega_R)}^2 
		\leq \Im(s)\Im\big(\left\langle \DtN u, \conj{u}\right\rangle_{S_R}\big).
	\end{align}
	\noindent\emph{Case a}: 
	If  $\Im s\ne0$, the sign~\eqref{eqn:DtN_Im_sign} of the imaginary part of $\DtN$
	and~\eqref{eqn:re_norm_u_esimate} result~in
	\begin{align*}
		\Re(s)\|\opA^{1/2}\nabla u\|_{L^2(\Omega_R)}^2 
		+\Re(s)|s|^2\|p^{1/2}u\|_{L^2(\Omega_R)}^2 
		< 0
		\qquad\text{or}\qquad u|_{S_R}\equiv 0.
	\end{align*}
	The left-hand side is non-negative as $\Re s \geq0$.
	Hence $u|_{S_R}\equiv 0$. Thus the extension by $u|_{B_R^+}\equiv 0$ solves~\eqref{eqn:FSP}
	by \Cref{lem:equivalence}. 
	The unique continuation principle%
	~\cite{AKS:UniqueContinuationTheorem1962,Wol:PropertyMeasuresRN1992}
	for globally Lipschitz $\opA\in W^{1,\infty}(\Omega;\mathbb{S}^\dim)$ 
	and the argumentation in~\cite{BCT:UniquenessTimeHarmonic2012} imply $u\equiv 0$.

	\medskip
	\noindent\emph{Case b}: For $\Im s=0$ and $\Re s>0$,~\eqref{eqn:re_norm_u_esimate} reveals
	$\|p^{1/2}u\|_{L^2(\Omega_R)}\leq0$ implying $u\equiv 0$ in $\Omega_R$ by~\ref{ass:C2}.
	This establishes uniqueness; further details are omitted.
\end{proof}
\Cref{lem:uniqueness} also holds for piecewise Lipschitz 
coefficients over certain countable sets $(\omega_j)_{j\in\N}$, 
see~\cite[Assumption~1.1]{BCT:UniquenessTimeHarmonic2012} 
and \cite[Prop.~2.13]{LRX:MoscoConvergenceHcurl2019} for details.
Lipschitz continuity is essentially optimal for uniqueness
in the indefinite case $\Re s=0$,
see~\cite{Fil:SecondOrderEllipticEquation2001} for an explicit counterexample with
$\alpha$-Hölder regular $\opA\in C^{0,\alpha}(\R^\dim;\mathbb{S}^\dim)$ 
for any $\alpha\in(0,1)$.

A consequence of the uniqueness by \Cref{ass:A_p} and Fredholm alternative arguments available 
for the truncated problem~\eqref{eqn:THP}
on the bounded domain $\Omega_R=\Omega\cap B_R$
is the existence of (unique) solutions.
The sesquilinear form $\opl:H^1(\Omega_R)\times H^1(\Omega_R)\to \R$ associated to~\eqref{eqn:THP}
is given for any $v,w\in H^1(\Omega_R)$ by
\begin{align}\label{eqn:l_def}
	\opl(v,w) \coloneqq\int_{\Omega_R}\left(\opA\nabla v\cdot \nabla \conj w 
			+ s^2 p\, v\,\conj w\right)\d x
		- \left\langle \DtN v, \conj w\right\rangle_{S_R}.
\end{align}
The dual space 
$\HdualGammaD(\Omega_R)
=(H^1_{\Gamma_{\mathrm{D}}}(\Omega_R))'$ 
is isomorphic to
$\{F\in(H^{1}_{\Gamma_{\mathrm{D}}}(\Omega))'\ 
:\ \mathrm{supp}(F)\subset \overline{\Omega_R}\}$.
The weak form of the truncated Helmholtz problem~\eqref{eqn:THP} 
for $F\in \HdualGammaD(\Omega_R)$ 
and $(\gD,\gN)\in H^{1/2}(\Gamma_{\mathrm{D}})\times H^{-1/2}(\Gamma_{\mathrm{N}})$
seeks $u\in %
H^1(\Omega_R)$ 
with $u|_{\Gamma_{\mathrm{D}}}=\gD$ and
\begin{align}\label{eqn:WTHP}
	\opl(u,v) = F(\conj v)+\left\langle \gN, \conj v\right\rangle_{\Gamma_{\mathrm{N}}}\qquad\text{for all
	}v\in \HGammaD(\Omega_R).
\end{align}
Let $\opL:\HGammaD(\Omega_R)\to \HdualGammaD(\Omega_R)$
denote the linear operator associated to $\opl$ by
\begin{align}\label{eqn:L_def}
	\opl(v,w) = \left\langle \opL v,\conj w\right\rangle_{\Omega_R}
		\qquad\text{for all }v,w\in \HGammaD(\Omega_R)
\end{align}
(in terms of the dual pairing
$\left\langle \bullet,\bullet\right\rangle_{\Omega_R}
=\left\langle \bullet,\bullet\right\rangle_%
{\HdualGammaD(\Omega_R)\times \HGammaD(\Omega_R)}$
from \Cref{sub:general_notation}).
A G\r arding inequality 
\cite{MS:ConvergenceAnalysisFinite2010,SW:WavenumberexplicitParametricHolomorphy2023} for $\opL$
implies the well-posedness of~\eqref{eqn:WTHP}, i.e.,
the continuity of the solution operator $\opN\coloneqq \opL^{-1}$.
Recall the weighted norm $\|\bullet\|_{H^1(\Omega_R),s}$ from
\eqref{eqn:H_norm_def} that induces the operator norm $\|\bullet\|_{\HdualGammaD(\Omega_R),s}$
for the dual space $\HdualGammaD(\Omega_R)$ by 
\begin{align}\label{eqn:H1_dualnorm_def}
	\|F\|_{\HdualGammaD(\Omega_R),s}
	\coloneqq\sup_{0\ne v\in H^{1}_{\mathrm{D}}(\Omega_R)}
	\frac{|F(\conj v)|}{\|v\|_{H^1(\Omega_R),s}}
		\qquad\text{for all }F\in \HdualGammaD(\Omega_R).
\end{align}
\begin{theorem}[existence and uniqueness
{\cite{MS:ConvergenceAnalysisFinite2010,SW:WavenumberexplicitParametricHolomorphy2023}}]
	\label{thm:L_N_well_def}
	Let \Cref{ass:A_p} be satisfied.
	The bounded operator 
	$\opL:\HGammaD(\Omega_R)\to \HdualGammaD(\Omega_R)$ 
	from~\eqref{eqn:DtN_non_pos} has a bounded inverse
	$\opN:\HdualGammaD(\Omega_R)\to \HGammaD(\Omega_R)$ with 
	\begin{align}\label{eqn:CN_def}
		\CN(s)
		\coloneqq \sup_{F\in \HdualGammaD(\Omega_R)}\,
		\frac{\|\opN F\|_{H^1(\Omega_R),s}}{\|F\|_{\HdualGammaD(\Omega_R),s}}
		<\infty.
	\end{align}
	In particular, 
	there exists a unique solution $u\in H^1(\Omega_R)$ 
	to~\eqref{eqn:WTHP} for any $F\in \HdualGammaD(\Omega_R)$.
\end{theorem}
\begin{proof}
	The properties of the $\DtN$ operator established for general 
	$s\in\Cstar_{\geq0}$ and $\dim\geq2$ in%
	~\cite[Thm.~3.3]{GS:DirichlettoNeumannOperatorHelmholtz2025}
	permit the application of the Fredholm alternative, following%
	~\cite[Sec.~3]{MS:ConvergenceAnalysisFinite2010} and%
	~\cite{GPS:HelmholtzEquationHeterogeneous2019,SW:WavenumberexplicitParametricHolomorphy2023}, as outlined below.
	The boundedness of $\DtN$ and the coefficients 
	(by~\ref{ass:C2}) implies the
	continuity of $\opL$.
	The G\r arding inequality
	\begin{align*}
		\Re(\opl(v,v))
		&\geq a_{\min}\|v\|_{H^1(\Omega_R)}^2 + 
		(\Re(s^2)p_{\min}-a_{\min})\|v\|_{L^2(\Omega_R)}^2
	\end{align*}
	holds by~\ref{ass:C2} and~\eqref{eqn:DtN_non_pos}.
	Since solutions to~\eqref{eqn:WTHP} are unique (by~\Cref{ass:A_p}), 
	the Fredholm alternative
	\cite[Thm.~2.34]{McL:StronglyEllipticSystems2000}
	verifies $\opL$ as a bounded linear bijection.
\end{proof}

\begin{remark}[coercivity of $\opl$ for $\Re s>0$]
	The proof of \Cref{thm:L_N_well_def} significantly simplifies in the case
	$\Re s>0$ with a coercive sesquilinear form $\opl$. Indeed, the coercivity 
	\begin{align*}
		\Re \left(\frac{\conj s}{|s|}\opl(v,v)\right)
		\geq\min\{a_{\min},p_{\min}\}\frac{\Re s}{|s|}\|v\|^2_{H^1(\Omega_R),s}
		\qquad\text{for all }v\in \HGammaD(\Omega_R)
	\end{align*}
	follows as in \cite[Lem.~3.2]{FHS:SkeletonIntegralEquations2024}
	and%
	~\cite{BH:FormulationVariationnelleEspacetemps1986,BS:IntegralEquationMethods2022}
	from $0\leq -\Re(\conj s\left\langle \DtN v,\conj v\right\rangle _{S_R})$.
	Hence the norm~\eqref{eqn:CN_def} of $\opN$ has the 
	upper bound (that degenerates as $\Re(s)\to 0$)
	\begin{align*}
		\CN(s)\leq \max\{a_{\min}^{-1},p_{\min}^{-1}\}\frac{|s|}{\Re(s)}.
	\end{align*}
\end{remark}
\begin{remark}[bounds on $\CN$ for $\Re s=0$]
	In the purely imaginary regime $\Re s=0$, the known upper bounds
	\cite{GPS:HelmholtzEquationHeterogeneous2019,SW:WavenumberexplicitParametricHolomorphy2023} for $\CN(s)$ depend
	polynomially on $\Im(s)$ for ``most frequencies''.
	However, there exist frequencies on the imaginary axis with 
	a super-algebraic growth of the operator norm $\CN(s)$
	\cite{PV:ResonancesRealAxis1999,GPS:HelmholtzEquationHeterogeneous2019}, i.e.,
	for all $m\in\mathbb{N}$ there is a constant $C_m>0$ and a sequence 
	$(s_n^m)_{n\in\mathbb{N}}\subset \operatorname{i}\R$ with 
	$C_m s_n^m\leq \CN(s_n^m)$ for all $n\in\mathbb{N}$.
\end{remark}

\subsection{The acoustic Helmholtz and solution operators}%
\label{sub:Newton_potential}
The remaining parts of this section analyse the acoustic Helmholtz and solution operators
$\opL$ and $\opN$.
It is known from%
~\cite{MS:ConvergenceAnalysisFinite2010,GS:DirichlettoNeumannOperatorHelmholtz2025}
that $\DtN$ coincides with its (linear) dual $\DtN'$.
Hence $\opL$ and its inverse $\opN$ are self-dual in the sense that
\begin{align}\label{eqn:L_N_dual_def}
	\left\langle \opL v,\conj w\right\rangle_{\Omega_R}
			= \left\langle v,\opLdual \conj{w}\right\rangle _{\Omega_R}
	\quad\text{and}\quad
	\left\langle \opN \varphi, \conj\psi\right\rangle_{\Omega_R} 
			= \left\langle \varphi, \opNdual \conj{\psi}\right\rangle_{\Omega_R}
\end{align}
holds for all $v,w\in \HGammaD(\Omega_R)$ and $\varphi,\psi\in \HdualGammaD(\Omega_R)$.
The restriction of $\opN$ onto more regular $L^2$ sources remains an 
isomorphism onto its image identified by the following theorem.
Define 
the vector space $V(\Omega_R,\opA,s)$ and
the exterior Neumann jump $\jumpNext\bullet$ by
\begin{align}\nonumber
	V(\Omega_R,\opA,s)
	&\coloneqq\{v\in \HGammaD(\Omega_R) \ :\ \Div(\opA\nabla v)\in L^2(\Omega_R),\;
	(\traceN[\Omega]v)|_{\Gamma_{\mathrm{N}}}=0,\;\jumpNext{v}=0 \},\\
	\jumpNext{\bullet}&\coloneqq \traceN - \DtN.\label{eqn:V_BR_def}%
\end{align}
Recall 
$\|\bullet\|_{H(\Omega_R,\Div),s}$ and $\|\bullet\|_{H^1(\Omega_R),s}$
from~\eqref{eqn:H_norm_def}--\eqref{eqn:Hdiv_norm_def}
and equip $V(\Omega_R,\opA,s)$ with 
\begin{align}\label{eqn:V_norm_def}
	\|v\|_{V(\Omega_R,\opA, s)}\coloneqq
	\sqrt{\|\opA\nabla v\|_{H(\Omega_R,\Div),s}^2 + \|v\|_{H^1(\Omega_R),s}^2
	}
	\qquad\text{for all }v\in V(\Omega_R,\opA,s).
\end{align}
We remark that $V(\Omega_R,\opA,s)$ and $V(\Omega_R,\opA,\conj s)$ do not coincide 
in general%
\footnote{
The identity ${\opSD[ \conj s]g}=\conj{\opSD[s] \conj g}$ 
by~\eqref{eqn:farfield_problem} 
implies $\DtN[\conj s] g=\conj{\DtN[s]\conj g}$
for all $g\in H^{1/2}(S_R)$.
}.
\begin{theorem}[{$L^2$} sources]\label{lem:Newton_potential_L2}
	The (not relabelled) restrictions
	\begin{align*}
		\opL&:V(\Omega_R,\opA,s)\to  L^2(\Omega_R)
		\quad\text{and its inverse}\quad
		\opN:L^2(\Omega_R)\to V(\Omega_R,\opA,s)%
	\end{align*}
	are well-defined bounded linear maps (in the norms of 
	$V(\Omega_R,\opA,s)$ %
	and $L^2(\Omega_R)$)
	with 
	\begin{align}\label{eqn:L_L2_repr}
		\opL v = -\Div(\opA \nabla v) + s^2p v \qquad\text{for all }v\in V(B_R,\opA,s).
\end{align}
\end{theorem}
\begin{proof}
	This proof considers the case $\Gamma_{\mathrm{D}}=\partial\Omega$ 
	(with $\Gamma_{\mathrm{N}}=\emptyset$)
	for a simpler exposition while the extension to 
	$\Gamma_{\mathrm{N}}\ne\emptyset$ is straightforward.
	Recall the exterior Neumann jump~\eqref{eqn:V_BR_def}.
	Given any $v\in \HGammaD(\Omega_R)= H^1_{\partial\Omega}(\Omega_R)$ 
	with $\Div(\opA\nabla v)\in L^2(\Omega_R)$,
	the definition of $\opl$ and $\opL$ 
	in~\eqref{eqn:l_def}--\eqref{eqn:L_def} plus 
	an integration by parts with arbitrary $w\in \HGammaD(\Omega_R)$ provide
	\begin{align}\label{eqn:range_L_L2}
		\left\langle \opL v, \conj w\right\rangle _{\Omega_R} = \opl(v,w) 
		=\int_{\Omega_R} (-\Div(\opA\nabla v) + s^2 p v) \conj w\d x +
	\left\langle \jumpNext v, \conj w\right\rangle_{S_R}. 
	\end{align}
	This and the vanishing jump $\jumpNext{v} \equiv (\traceN -\DtN)(v) = 0$
	and $\Div(\opA \nabla v)\in L^2(\Omega_R)$
	for any $v\in V(\Omega_R,\opA,s)$
	verify~\eqref{eqn:L_L2_repr} for $\opL v\in L^2(\Omega_R)$.
	To prove surjectivity, let $f\in L^2(\Omega_R)$ be arbitrary and set 
	$v\coloneqq \opN f\in \HGammaD(\Omega_R)$. 
	An integration by parts with an arbitrary $w\in C^\infty_0(\Omega_R)$
	and $\opL v = \opL\opN f = f\in L^2(\Omega_R)$
	reveal
	\begin{align*}
		\int_{\Omega_R}\opA \nabla v\cdot \nabla \conj w\d x 
		= \left\langle \opL v, \conj{w}\right\rangle_{\Omega_R} 
		- \int_{\Omega_R}s^2p\,v\, \conj w\d x = \int_{\Omega_R} (f - s^2 p v)
		\conj w \d x.
	\end{align*}
	Hence, the weak divergence
	$-\Div(\opA\nabla v) = f - s^2 p v\in L^2(\Omega_R)$ is square-integrable.
	Moreover, $\opL v = f = -\Div(\opA\nabla v) + s^2 p v \in L^2(\Omega_R)$ 
	and~\eqref{eqn:range_L_L2}
	verify 
	\begin{align*}
		\left\langle(\traceN -\DtN)v, \conj w \right\rangle_{S_R} = 0 
		\qquad\text{for all }w\in 
		\HGammaD(\Omega_R).
	\end{align*}
	Since the jump $\jumpNext{v}=0$ must vanish by the fundamental theorem 
	of the calculus of variations, this shows $v\in V(\Omega_R,\opA,s)$ and 
	implies the surjectivity
	$\opL:V(\Omega_R,\opA,s)\to L^2(\Omega_R)$ so that the Newton potential 
	$\opN:L^2(\Omega_R)\to V(\Omega_R,\opA,s)$ is well defined.
	
	The boundedness of $\opL:V(\Omega_R,\opA,s)\to L^2(\Omega_R)$ follows immediately 
	from~\eqref{eqn:L_L2_repr} and the definition~\eqref{eqn:V_norm_def}
	of the norm in $V(\Omega_R,\opA,s)\subset \HGammaD(\Omega_R,\opA)$.
	By the open mapping theorem, the inverse $\opN:L^2(\Omega_R)\to V(\Omega_R,\opA,s)$ 
	is bounded as well. 
\end{proof}

Define the (weighted)
operator norm of $\opN$ and $\opNdual$ in $L(L^2(\Omega_R); V(\Omega_R,\opA,s))$~as
\begin{align}\label{eqn:CN_L2_def}
	\tCN(s)
	\coloneqq \sup_{0\ne f\in L^2(\Omega_R)}\,
	\frac{\|\opN f\|_{V(\Omega_R,\opA,s)}}{|s|^{-1}\|f\|_{L^2(\Omega_R)}}.
\end{align}
The scaling in the wavenumber $|s|$ in~\eqref{eqn:CN_L2_def} matches that 
of $\CN(s)$ from~\eqref{eqn:CN_def}.
\begin{lemma}[bound on $\tCN(s)$]\label{lem:tCN_bound}
	It holds %
	$\tCN(s)\leq \sqrt{2+(2p_{\max}^2+1)\CN^2(s)}$.
\end{lemma}
\begin{proof}
	Let $f\in L^2(\Omega_R)$ be arbitrary and set $v\coloneqq\opN f\in V(\Omega_R,\opA,s)$.
	A triangle inequality and $\opL v=f$ with~\eqref{eqn:L_L2_repr} reveals
	with~\ref{ass:C2} that
	\begin{align*}
		\|\Div(\opA\nabla v)\|_{L^2(\Omega_R)}
		\leq \|f\|_{L^2(\Omega_R)}+|s|^2\|pv\|_{L^2(\Omega_R)}
		\leq \|f\|_{L^2(\Omega_R)}+|s|^2p_{\max}\|v\|_{L^2(\Omega_R)}.
	\end{align*}
	Hence, the definition of 
	$\|\bullet\|_{H^1(\Omega_R),s}$ and
	$\|\bullet\|_{V(\Omega_R,\opA,s)}$ 
	in~\eqref{eqn:H_norm_def} and~\eqref{eqn:V_norm_def}
	result in
	\begin{align*}
		\|v\|_{V(\Omega_R,\opA,s)}^2
		&\leq (2p_{\max}^2+1)\|v\|_{H^1(\Omega_R),s}^2
		+ 2|s|^{-2}\|f\|_{L^2(\Omega_R)}^2\\
		&\leq (2p_{\max}^2+1)\CN^2(s)\|f\|_{\HdualGammaD(\Omega_R),s}^2
		+ 2|s|^{-2}\|f\|_{L^2(\Omega_R)}^2
	\end{align*}
	with the operator norm $\CN(s)$ of $\opN\in L(\HdualGammaD(\Omega_R);\HGammaD(\Omega_R))$ from~\eqref{eqn:CN_def}
	in the last step.
	This and $|s|\|f\|_{\HdualGammaD(\Omega_R),s}\leq\|f\|_{L^2(\Omega_R)}$ by definition conclude the proof.
\end{proof}
\begin{remark}[comparison with~\cite{FHS:SkeletonIntegralEquations2024}]\label{rem:B_infty}
	The reformulation of the exterior Helmholtz problem~\eqref{eqn:FSP}
	on the truncated domain $\Omega_R=\Omega\cap B_R$ appears necessary for purely imaginary 
	Helmholtz problems with $\Re s=0$ and enables a unified analysis for general
	wavenumbers $s\in \Cstar_{\geq0}$.

	For wavenumbers with positive real part $\Re s>0$, all solutions to the 
	full Helmholtz problem~\eqref{eqn:FSP} 
	satisfy the integrability $u\in H^1(\Omega)$ 
	over the whole computational domain $\Omega$
	and the truncation of the computational domain is not necessary for the analysis.
	Indeed, for $\Re s>0$, the analysis in this paper applies also to $R=\infty$
	with the conventions $\Omega_\infty=\Omega,B_\infty=\R^n$, and $S_\infty=\emptyset$
	such that the truncated Helmholtz problem~\eqref{eqn:THP} 
	coincides with~\eqref{eqn:FSP}.
	In this case, the conditions \ref{ass:C3}--\ref{ass:C4} are redundant and 
	the well-posedness follows from the coercivity~\cite[Lem.~3.2]{FHS:SkeletonIntegralEquations2024}
	of the associated sesquilinear
	form~\eqref{eqn:l_def}, while the results in the subsequent sections recover 
	and overcome the limitations in~\cite{FHS:SkeletonIntegralEquations2024}.
\end{remark}
\section{Potential operators for interface problems}%
\label{sec:Potential operators for interface problems}
The solution operator from~\Cref{sec:potential_operators}
enables a variational definition of single, double, and boundary layer potentials
for the Helmholtz operator with varying coefficients.

\subsection{The transmission problem for a single interface}%
\label{sub:interface problem}
The compact interface $\Interface\coloneqq \partial \IGebiet$
is the boundary of either some (in particular connected) exterior Lipschitz 
domain $\IGebiet\subset\R^n$ or
some bounded (possibly multiply connected) Lip\-schitz set $\IGebiet\subset\R^n$.
Throughout this section, the computational domain is the full space $\Omega=\R^n$.
The wavenumber $s\in \Cstar_{\geq0}$ 
and the coefficients $\opA\in L^\infty(\R^n;\mathbb{S}^n)$ and 
$p\in L^\infty(\R^n)$ satisfy~\ref{ass:C1}--\ref{ass:C3} for a 
sufficiently large ball $B_R\subset\R^n$ that contains $\Interface \subset B_R$.
We require the analogue of \Cref{ass:A_p} in the current setting.
\begin{assumption}\label{ass:A_p_new}
	For any $F\in H^{-1}_{\mathrm{D, comp}}(\R^n)$ with~\ref{ass:C4},
	there exists at most one
	solution to the truncated problem~\eqref{eqn:THP} on $\Omega_R=B_R$ 
	(with $\Gamma_{\mathrm{D}}=\emptyset=\Gamma_{\mathrm{N}}$).
\end{assumption}
The transmission problem on $\Interface$ 
seeks a weak
solution $u\in H^1_{\mathrm{loc}}(\R^\dim\setminus\Interface)$ to
\begin{equation}\label{eqn:IP}
	\begin{aligned}
		-\Div(\opA \nabla u) + s^2 p\, u 
			&=0&&\text{in }\R^\dim\setminus\Interface,\\
		\span\span u\text{ satisfies }\eqref{radicond}\span
	\end{aligned}
\end{equation}
with prescribed %
jumps $\jumpD[\IGebiet]u=\gD\in H^{1/2}(\Interface)$ and 
$\jumpN[\IGebiet]u=\gN\in H^{-1/2}(\Interface)$
across $\Interface$.
Solutions to~\eqref{sub:interface problem} are characterised
in the exterior domain $B_R^+=\R^n\setminus\overline{B_R}$
by~\eqref{eqn:farfield_problem}.
Hence their restrictions to $B_R$ lie in the space
$V(B_R\setminus\Interface,\opA,s)$ 
defined in analogy to~\eqref{eqn:V_BR_def} by
\begin{align}\label{eqn:V_BR_interface_def}
	V(B_R\setminus\Interface,\opA,s)
	&\coloneqq\{v\in H^1(B_R\setminus\Interface,\opA) \ :\ \jumpNext{v}=0 \}
\end{align}
with the norm $\|\bullet\|_{V(B_R\setminus\Interface,\opA,s)}$
as in~\eqref{eqn:V_norm_def} (for $\Omega_R$ replaced by $B_R\setminus\Interface$).
The equivalent formulation of~\eqref{eqn:IP} (in the sense of~\Cref{lem:equivalence}) seeks a
solution $u\in V(B_R\setminus\Interface,\opA,s)$ to%
\begin{subequations}\label{eqn:TIP}
	\begin{align}\label{eqn:TIPa}
		-\Div(\opA \nabla u) + s^2 p\, u 
			&=0&&\text{in }B_R\setminus\Interface,\\
\jumpD[\IGebiet]u = \gD
		\quad\text{and}\quad
		\jumpN[\IGebiet]u &= \gN&&\text{on }\Interface\label{eqn:TIPb}
	\end{align}
\end{subequations}
for given Dirichlet data $\gD\in H^{1/2}(\Interface)$ and Neumann data $\gN\in H^{-1/2}(\Interface)$.
(The condition $\partial_r u = \DtN u$ on $S_R=\partial B_R$ is implied by $u\in V(B_R\setminus\Interface,\opA,s)$).
This section introduces and analyses an integral formulation of~\eqref{eqn:TIP} 
based on a novel variational definition of the single layer 
potential $\opS$ and 
the double layer potential $\opD$
extending the approach for the coercive case (that is $\Re s>0$)
in~\cite{FHS:SkeletonIntegralEquations2024} to purely imaginary wavenumbers.

To provide a sharp wavenumber-explicit stability analysis, 
the remaining parts of this subsection discuss weighted trace
 norms introduced and analysed in~\cite{Gra:OptimalTraceNorms2025}.
Let 
$\IGebietR\coloneqq \IGebiet\cap B_R$
denote the intersection of $\IGebiet$ with $B_R$.
The trace space $H^{1/2}(\Interface)$ is naturally equipped with the minimal extension norm
\begin{align}\label{eqn:H12_norm_def}
	\NormHD[\Interface]{s}g
	\coloneqq\inf_{\substack{v\in H^1(\IGebietR)\\[0.2em] \traceD[\IGebiet]v = g}}
		\|v\|_{H^1(\IGebietR),s}
	\qquad\text{for all }g\in H^{1/2}(\Interface).
\end{align}
This trace norm arises naturally from the identification of $H^{1/2}(\Interface)$ with the quotient space%
\footnote{
	The identification with the quotient space
	can further be utilised to define abstract trace spaces
	in a generalised setting with non-Lipschitz $\IGebietR$ 
	that does not admit classical traces%
	~\cite{CH:IntegralEquationsMultiScreens2013,HPS:TracesHilbertComplexes2023}.
}
$H^1(\IGebietR)/H^1_{\Gamma}(\IGebietR)$ 
equipped with the energy norm~\eqref{eqn:H_norm_def}.
An intrinsic characterisation of~\eqref{eqn:H12_norm_def} in terms of a weighted Sobolev-Slobodeckij-type norm is
provided in~\cite[Sec.~3]{Gra:OptimalTraceNorms2025}. 
The dual space $H^{-1/2}(\Interface)=(H^{1/2}(\Interface))'$ 
is equipped with the operator norm
\begin{align}\label{eqn:H12_dual_norm_def}
	\NormHN[\Interface]{s}{{h}}
	&\coloneqq
	\sup_{\substack{0\ne {g}\in H^{1/2}(\Interface)}}
	\frac{\left|\left\langle {h},{g}\right\rangle_{\Interface}\right|}
	{\NormHD[\Interface]{s}{{g}}}
	\qquad\text{for all }{h}\in H^{-1/2}(\Interface).
\end{align}
For $s=1$,~\eqref{eqn:H12_norm_def}--\eqref{eqn:H12_dual_norm_def} are classical trace norms and equivalent, e.g., to
Sobolev-Slobodeckij or interpolation norms \cite{LM:NonhomogeneousBoundaryValue1972}.
Their scaling in the weight $s$
is identified in~\cite[Lem.~4.1]{Gra:OptimalTraceNorms2025}
for any $g\in H^{1/2}(\Gamma)$ and $h\in H^{-1/2}(\Gamma)$ as%
\begin{equation*} 
	\begin{aligned}
		\min\{1,|s|\}\NormHD[\Gamma]{1}{g}
		&\leq
		\NormHD[\Gamma]{s}{g}
		\leq \max\{1,\min\{|s|,\Const{sc}|s|^{1/2}\}\}\NormHD[\Gamma]{1}{g},\\
		\min\{1,|s|\}\NormHN[\Gamma]{s}{h}
		&\leq
		\NormHN[\Gamma]{1}{h}
		\leq\max\{1,\min\{|s|,\Const{sc}|s|^{1/2}\}\}\NormHN[\Gamma]{s}{h}
	\end{aligned}
\end{equation*}
with some universal constant $\Const{sc}>0$.
The following result recalls the $s$-explicit trace inequality from%
~\cite{Gra:OptimalTraceNorms2025} in terms of a lower bound on the wavenumber modulus%
\begin{align}\label{eqn:s_low_def}
	\slow\coloneqq\min\{1,|s|\}\leq 1
	\qquad\text{and}\qquad
	\sLow\coloneqq\slow^{-1}=\max\{1,|s|^{-1}\}\geq1.
\end{align}%
The properties of the trace 
norms~\eqref{eqn:H12_norm_def}--\eqref{eqn:H12_dual_norm_def}
depend on the geometry of the extension set $\IGebietR$
and its boundary $\partial\IGebietR\subset\Gamma\cup S_R$.
By assumption on $\IGebiet$, either $\IGebiet\subset B_R$ or
its complement $\R^\dim\setminus\overline \IGebiet\subset B_R$ 
is bounded.
Denote this bounded set by
$\IGebietInt\subset B_R$ (with $\Interface=\partial\IGebietInt=\partial \IGebiet$).
\begin{lemma}[$s$-explicit trace estimate]\label{lem:trace_inequality}%
	There exist constants $\Const{tr,D},\Const{tr,N}>0$ independent of $s$ and 
	exclusively depend on $\IGebietR$ and $\Gamma$ with
	\begin{align}\label{eqn:Ctr_D}
		\|\traceD[\IGebiet] v\|_{H^{1/2}(\Interface),s}
		&\leq \|v\|_{H^1(\IGebietR),s}
		&&\text{for all }v\in H^1(\IGebietR),\\\label{eqn:Ctr_N}
		\NormHN[\Interface]{s}{\traceNu[\IGebietInt]{\mathbf{q}}}
		&\leq \Const{tr,N}\|\mathbf{q}\|_{H(\IGebietInt,\Div),s}
		&&\text{for all }\mathbf{q}\in H(\IGebietInt,\Div).
	\end{align}
	Moreover, any $v\in H^1(B_R\setminus\overline{\IGebietR})$ and 
	$\mathbf{q}\in H(B_R\setminus\overline{\IGebietInt},\Div)$
	satisfy
	\begin{align}\label{eqn:Ctr_Db}
		\|\traceDext[\IGebiet] v\|_{H^{1/2}(\Interface),s}
		&\leq \Const{tr,D} \|v\|_{H^1(B_R\setminus\overline{\Omega_R}),\max\{1,|s|\}}
		\leq \Const{tr,D}\sLow\, \|v\|_{H^1(B_R\setminus\overline{\Omega_R}),s}, \\
		\NormHN[\Interface]{s}{\traceNuext[\IGebietInt]{\mathbf{q}}}
		&\leq \Const{tr,N}\sLow\, \|\mathbf{q}\|_{H(B_R\setminus\overline{\IGebietInt},\Div),s}.\label{eqn:Ctr_Nb}
	\end{align}
	If $\IGebietR=\IGebietInt$,~\eqref{eqn:Ctr_N} holds for $\Const{tr,N}$ 
	replaced by $1$.
\end{lemma}

\begin{proof}
	This follows from a straightforward distinction between the two cases
	$\IGebietR=\IGebietInt$
	with $\partial\IGebietR=\Gamma$
	and $\IGebietR\ne\IGebietInt$
	(with $\partial\IGebietInt=\Gamma$)
	from~\cite[Thm.~4.4]{Gra:OptimalTraceNorms2025} in the current setting; further details are omitted.
\end{proof}

\subsection{Single layer potential}%
\label{sub:Single_layer_potential}
The solution operator $\opN$ from \Cref{thm:L_N_well_def} for $\Omega_R=B_R$
in the current setting is also called \emph{acoustic Newton potential} in the following.
Let $\traceD[\IGebiet]':H^{-1/2}(\Interface)\to \Hdual(B_R)$ 
denote the dual operator of the Dirichlet trace map $\traceD[\IGebiet]$
from $H^1(B_R)$ onto the interface $\Gamma=\partial\IGebiet$.
In analogy to the classical definition in~\cite[p.~202]{McL:StronglyEllipticSystems2000}
(and also~\cite[Def.~3.1.5]{SS:BoundaryElementMethods2011}),
the single layer potential is defined in the present situation as the composition
\begin{align}\label{eqn:S_def}
	\opS\coloneqq \opN\traceD[\IGebiet]':H^{-1/2}(\Interface)\to H^1(B_R).
\end{align}
The variational formulation of $\opS$ from~\eqref{eqn:S_def} reads
\begin{align}\label{eqn:S_def_var}
	\boxed{
	\opl(\opS \tracevarb, v) 
	= \left\langle \tracevarb, \traceD[\IGebiet] \conj{v}\right\rangle_{\Interface}
		\qquad\text{for all }\tracevarb\in H^{-1/2}(\Interface), v\in H^1(B_R)%
	}
\end{align}%
and has been previously used to define and analyse generalised single layer potentials,
see, e.g.,%
~\cite[Sec.~4--5]{Bar:LayerPotentialsGeneral2017}
for an abstract setting
and~\cite[Def.~3.6]{FHS:SkeletonIntegralEquations2024} for the definite case ($\Re s>0$).
\begin{theorem}[single layer potential]\label{lem:S_property}
	The operator $\opS$ from~\eqref{eqn:S_def} maps $H^{-1/2}(\Interface)$ 
	boundedly into $H^1(B_R)\cap V(B_R\setminus\Interface,\opA,s)$ and 
	is uniquely defined by~\eqref{eqn:S_def_var}.
	Any $g\in H^{-1/2}(\Interface)$ and $u\coloneqq \opS g$ satisfy for
	$\Const{SL}\coloneqq(1+\max\{a_{\max},p_{\max}\}^2)^{1/2}$ that
	\begin{enumerate}[label=(\roman*)]
		\item $\displaystyle
			\Const{SL}^{-1}
			\|u\|_{V(B_R\setminus\Interface,\opA, s)}\leq \|u\|_{H^1(B_R),s}\leq
		\CN(s)\|g\|_{H^{-1/2}(\Interface),s}$,
		\item $-\Div(\opA \nabla u) + s^2 p\, u = 0$ in $B_R\setminus\Interface$ and
		\item $\jumpD[\IGebiet]{u} = 0
		\quad\text{and}\quad
		\jumpN[\IGebiet]{u} = \tracevarb$.
	\end{enumerate}
\end{theorem}
\begin{proof}
	The equivalence of~\eqref{eqn:S_def}--\eqref{eqn:S_def_var} is clear.
	Consider $u\coloneqq\opS \tracevarb$ for any $\tracevarb\in H^{-1/2}(\Interface)$ 
	and let $\oldvarphi\in C^\infty_0(\overline{B_R}\setminus\Interface)$ be arbitrary.
	Since $\traceD[\IGebiet]\conj \oldvarphi =0$,
	an integration by parts with~\eqref{eqn:l_def},~\eqref{eqn:V_BR_def},
	and \eqref{eqn:S_def} verify
	\begin{align*}
		0=\opl(u,\oldvarphi) = 
		\int_{B_R}\big(-\Div(\opA \nabla u)\conj\oldvarphi 
		+ s^2p\, u\,\conj\oldvarphi\big)\d x + 
		\left\langle\jumpNext{u}, \conj
		\oldvarphi \right\rangle_{S_R}.
	\end{align*}
	As in the proof of \Cref{thm:L_N_well_def}, this reveals $\jumpNext{u} = 0$ and
	$-\Div(\opA\nabla u) = -s^2p\, u\in L^2(B_R\setminus\Interface)$ implying~$(ii)$.
	Hence, $u\in V(B_R\setminus\Interface,\opA,s)$ and (by~\ref{ass:C2}) 
	$\|\Div(\opA \nabla u)\|_{L^2(B_R\setminus\Interface)} \leq |s|^2
	p_{\max}\|u\|_{L^2(B_R)}$. %
	This,
	the definition~\eqref{eqn:V_norm_def} of the norm, and
	\ref{ass:C2} reveal
	\begin{align}\nonumber
		\|u\|_{V(B_R\setminus\Interface,\opA,s)}^2
		&\leq (1+a_{\max}^2)\|\nabla u\|_{L^2(B_R\setminus\Interface)}^2
			+ (1+p_{\max}^2)\|u\|_{L^2(B_R)}^2\\
		&\leq (1+\max\{a_{\max},p_{\max}\}^2)\|u\|_{H^1(B_R),s}^2.
		\label{eqn:SL_u_VBR_bound}
	\end{align}
	Since %
	$\|\traceD[\IGebiet]v\|_{H^{1/2}(\Interface),s}
	\leq \|v\|_{H^1(\IGebietR),s}
	\leq \|v\|_{H^1(B_R),s}$
	for all $v\in H^1(B_R)$ by \Cref{lem:trace_inequality},
	the operator norm of $\traceD[\IGebiet]':H^{-1/2}(\Interface)\to \Hdual(B_R)$
	is bounded above by $1$.
	Hence~\eqref{eqn:SL_u_VBR_bound},
	the characterisation $\opS=\opN\traceD[\IGebiet]'$ %
	and the definition of $\CN(s)$ in~\eqref{eqn:CN_def} reveal
	$(i)$.

	Similarly, the integration by parts formula for any $\oldvarphi\in C^\infty_0(B_R)$ 
	and~$(ii)$ provide
	\begin{align*}
		\left\langle \jumpN[\Interface]{u},\conj \oldvarphi\right\rangle_{\Interface} 
			= \left\langle \tracevarb, \conj{\oldvarphi}\right\rangle_{\Interface}.
	\end{align*}
	This and $\jumpD[\IGebiet]{u}=0$ from $u\in H^1(B_R)$ 
	verify~$(iii)$ and conclude the proof.
\end{proof}

\subsection{Double layer potential}%
\label{sub:Double_layer_potential}
The second ingredient for interface problems is the double layer potential 
that provides a solution of the homogeneous Helmholtz problem~\eqref{eqn:FSP} with 
prescribed Dirichlet jumps across $\Interface$. 
The double layer potential is classically defined~\cite[p.~202]{McL:StronglyEllipticSystems2000} 
as the composition of the full-space Newton potential and the dual Neumann trace
in analogy to~\eqref{eqn:S_def}.

Since the dual operator $\traceN[\IGebiet]'$ of the 
Neumann map
\begin{align*}
	\traceN[\IGebiet]:
	H^1(B_R,\opA)\to H^{-1/2}(\Interface)
\end{align*}
maps $H^{1/2}(\Gamma)$ into the dual space $(H^1(B_R,\opA))'$ 
which is strictly larger than
the domain of definition for $\opN$ from \Cref{thm:L_N_well_def},
that composition relies on an appropriate extension of $\opN$.
Recall the restriction of $\opN$ to $L^2(B_R)$ from \Cref{lem:Newton_potential_L2}.
Its dual operator $\opNddual:(V(B_R,\opA,s))'\to L^2(B_R)$ is given by
\begin{align}\label{eqn:Nddual_def}
	\left\langle \opNddual F, \conj{f}\right\rangle_{B_R} 
	\coloneqq \left\langle F, \opNdual \conj{f}\right\rangle _{B_R}
	\qquad\text{for all }F\in (V(B_R,\opA,s))', f\in L^2(B_R).
\end{align}
By the self-duality~\eqref{eqn:L_N_dual_def} of $\opN$ and the density $L^2(B_R)\subset \Hdual(B_R)$, 
this operator is indeed an extension and we write $\opN\coloneqq\opNddual$ 
in the following.
This and $(H^1(B_R,\opA))'\subset (V(B_R,\opA,s))'$
justifies the definition of the double layer potential as
\begin{align}\label{eqn:D_def}
	\opD\coloneqq \opN \traceN[\IGebiet]':H^{1/2}(\Interface)\to L^2(B_R).
\end{align}
Observe that $\opL[\conj s]:V(B_R,\opA,\conj s)\to L^2(B_R)$ is surjective
with $\conj{\opL[\conj s]v} = \opL\conj v$ for all $v\in V(B_R,\opA,\conj s)$
by \Cref{lem:Newton_potential_L2}.
Hence~\eqref{eqn:Nddual_def} and $\opN\opL =\id$ reveal an equivalent
variational characterisation of~\eqref{eqn:D_def} as
\begin{align}\label{eqn:D_ultraweak_def}
	\boxed{
	\left\langle \opD\tracevara, \opL\conj v\right\rangle_{B_R}
	= \left\langle \tracevara, \traceN[\IGebiet] \conj
	v\right\rangle_{\Interface}\qquad\text{for all }\tracevara\in H^{1/2}(\Interface), v\in V(B_R,\opA,\conj s).
}
\end{align}%
This generalises the variational definition in%
~\cite[Eqn.~(3.22)]{FHS:SkeletonIntegralEquations2024} for $\Re s>0$ and $R=\infty$.

The analysis of the double layer potential~\eqref{eqn:D_def} 
extends~\cite{FHS:SkeletonIntegralEquations2024}
based on a mixed reformulation of \eqref{eqn:D_ultraweak_def} 
with a separate variable for the weak gradient in $H(B_R,\Div)$:
Given any $g\in H^{1/2}(\Interface)$, 
seek $(\mathbf{p},\nomathfrak{u})\in M\coloneqq H(B_R,\Div)\times L^2(B_R)$ with
\begin{equation}\label{eqn:D_mixed_def}
	\begin{aligned}
		-\left\langle \opA^{-1}\mathbf{p},\conj{\mathbf{q}}\right\rangle_{B_R}
		+\big\langle \traceNu[B_R]\mathbf{p}, \conj{\NtD[\conj s]\traceNu[B_R]\mathbf{q}}
			\big\rangle_{S_R}
		- \left\langle \nomathfrak{u},\Div\conj{\mathbf{q}}\right\rangle_{B_R}
		&= \left\langle\tracevara, \conj{\traceNu[\IGebiet]\mathbf{q}}\right\rangle_{\Interface},\\
		-\left\langle \Div \mathbf{p}, \conj {\nomathfrak{v}}\right\rangle_{B_R} 
		+ \left\langle s^2p\,\nomathfrak{u},\conj{{\nomathfrak{v}}}\right\rangle
		_{B_R} 
		&= 0
	\end{aligned}
\end{equation}
for any $(\mathbf{q}, \nomathfrak{v})\in M$.
The weighted norm in $M$ is given by
\begin{align}\label{eqn:M_norm_def}
	\|(\mathbf{q},v)\|_{M,s}
	&\coloneqq \sqrt{\|\mathbf{q}\|_{H(B_R,\Div),s}^2+|s|^2\|v\|_{L^2(B_R)}^2}
	\qquad\text{for all }(\mathbf{q},v)\in M.
\end{align}%
The following lemma states the equivalence of~\eqref{eqn:D_ultraweak_def} 
and~\eqref{eqn:D_mixed_def}
as part \emph{(ii)} and extends the corresponding result~\cite[Lem.~3.10]{FHS:SkeletonIntegralEquations2024}.
\begin{lemma}[{mixed formulation}]
	\label{lem:well-posedness of mixed formulation}
	Given any $\tracevara\in H^{1/2}(\Interface)$,
	the mixed problem~\eqref{eqn:D_mixed_def} admits a unique solution 
	$(\mathbf{p}, u)\in M$.
	This unique solution satisfies
	\begin{enumerate}[label=(\roman*)]
		\item $\|(\mathbf{p},u)\|_{M,s}\leq \Const{DL}\left(1+\CN(s)\right)
			\|g\|_{H^{1/2}(\Interface),s}$,
		\item $u$ solves~\eqref{eqn:D_ultraweak_def} in place of $\opD\tracevara$,
		\item $u\in V(B_R\setminus\Interface,\opA,s)$ and $\mathbf{p}|_{B_R\setminus\Interface}=\opA\nabla u|_{B_R\setminus\Interface}
			\in H(B_R\setminus\Interface,\Div)$,
		\item $\jumpD[\IGebiet] u= -\tracevara$
			and $\jumpN[\IGebiet] u = 0$.
	\end{enumerate}
	The constant $\Const{DL}>0$ 
	exclusively depends on $a_{\max},p_{\max}$, and $p_{\min}$.
\end{lemma}
\begin{proof}
	The proof in two steps starts with the analysis of the mixed 
	formulation~\eqref{eqn:D_mixed_def}.

	\medskip
	\noindent\emph{Step 1} (well-posedness of~\eqref{eqn:D_mixed_def}):
	The sesquilinear form $b:M \times M\to \C$ 
	corresponding to~\eqref{eqn:D_mixed_def} 
	is given for any $\mathbf{p},\mathbf{q}\in H(B_R,\Div)$ and
	${u},{v}\in L^2(B_R)$ by
	\begin{equation}\label{eqn:b_def}
		\begin{aligned}
			b( (\mathbf{p},\nomathfrak{u}),(\mathbf{q},\nomathfrak{v}))
			\coloneqq 
			&- \left\langle \opA^{-1}\mathbf{p},\conj{\mathbf{q}}\right\rangle_{B_R}
			+ \big\langle \traceNu[B_R]\mathbf{p}, \conj{\NtD[\conj s]\traceNu[B_R]\mathbf{q}}
			\big\rangle_{S_R} 
			-\left\langle \nomathfrak{u},\Div\conj{\mathbf{q}}\right\rangle_{B_R}\\
			&-\left\langle \Div\mathbf{p},\conj{\nomathfrak{v}}\right\rangle_{B_R}
			+\left\langle s^2 p \nomathfrak{u},\conj{\nomathfrak{v}}\right\rangle_{B_R}.
		\end{aligned}
	\end{equation}
	To prove an inf-sup condition for $b(\bullet,\bullet)$,
	let $(\mathbf{p}, u)\in M$ be arbitrary and consider
	${w}\coloneqq \opNdual u\in V(B_R,\opA,\conj s)$.
	The definition of the
	norm $\|\bullet\|_{V(B_R,\opA,s)}$ in~\eqref{eqn:V_norm_def} and~\eqref{eqn:CN_L2_def} reveal
	\begin{align}\label{eqn:b_infsup_U_bound}
		\|\opA \nabla {w}\|_{H(B_R,\Div),s}^2+|s|^2\|{w}\|_{L^2(B_R)}^2
		\leq \|{w}\|_{V(B_R,\opA,s)}^2\leq \tCN(s)^2\,|s|^{-2}\,\|u\|_{L^2(B_R)}^2.
	\end{align}
	Since $\mathrm{NtD}$ is the inverse of $\mathrm{DtN}$, 
	$(\traceSD - \NtD[\conj s]) \varphi=-\NtD[\conj s]\jumpNextdual{\varphi}= 0$
	holds for all $\varphi\in V(B_R,\opA,\conj s)$ and
	the integration by parts formula verifies
	\begin{equation*}
		\begin{aligned}
			-\left\langle \Div \mathbf{p}, \conj {w}\right\rangle_{B_R}
			&= - \left\langle \traceNu[B_R]\mathbf{p}, \conj{{w}}\right\rangle_{S_R}
			+ \left\langle \mathbf{p}, \nabla\conj{{w}}\right\rangle_{B_R}\\
			&= - \big\langle \traceNu[B_R]\mathbf{p}, 
			\conj{\NtD[\conj s] {w}}\big\rangle_{S_R}
			+ \big\langle \opA^{-1}\mathbf{p}, \opA\nabla\conj{{w}}\big\rangle_{B_R}.
		\end{aligned}
	\end{equation*}
	(Recall from \Cref{sub:The full space Helmholtz problem} 
	that we abbreviate $\NtD[\conj s]{w}\coloneqq \NtD[\conj s]\traceN {w}$.)
	This, 
	\eqref{eqn:b_def} for
	$\mathbf{q}\coloneqq\opA\nabla {w}\in H(B_R,\Div)$, and $v={w}\in L^2(B_R)$
	reveal with~\eqref{eqn:L_L2_repr} that
	\begin{align}\label{eqn:jv_ibp}
		b( (\mathbf{p},\nomathfrak{u}), (\opA\nabla {w}, {w}))
		= \left\langle \nomathfrak{u}, \opL \conj {w}\right\rangle_{B_R}
		= \|u\|_{L^2(B_R)}^2%
	\end{align}
	with
	$\opL\conj{{w}}=\conj{\opLdual\opNdual u}=\conj u$ by \eqref{eqn:L_L2_repr}
	in the last step.
	Elementary algebra reveals
	\begin{align*}
		b\left( (\mathbf{p},\nomathfrak{u}), (-\mathbf{p}, u)\right)
			 &= \|\opA^{-1/2}\mathbf{p}\|_{L^2(B_R)}^2 
			 - \big\langle \traceNu[B_R]\mathbf{p}, \conj{\NtD[\conj s]\traceNu[B_R]\mathbf{p}}\big\rangle_{S_R}\\
			 &\quad+2i\Im\big(\left\langle u,\Div{\conj{\mathbf{p}}}\right\rangle_{B_R}\big)
			 +s^2\|p^{1/2}u\|_{L^2(B_R)}^2,\\
		 b((\mathbf{p},\nomathfrak{u}), (0, s^2 u+p^{-1}\Div \mathbf{p}))
			 &= -\|p^{-1/2}\Div\mathbf{p}\|_{L^2(B_R)}^2 \\
			 &\quad+2i\Im\big(\left\langle s^2u,
				\Div\conj{\mathbf{p}}\right\rangle_{B_R}\big)
				 + |s|^4 \|p^{1/2}u\|_{L^2(B_R)}^2.
	\end{align*}
	The real part of $\DtN$ is non-positive~\eqref{eqn:DtN_non_pos}
	by~\cite{Ned:AcousticElectromagneticEquations2001,MS:ConvergenceAnalysisFinite2010,GS:DirichlettoNeumannOperatorHelmholtz2025}.
	Hence,
	\begin{align*}
		0\leq-\Re\big(\left\langle \NtD \tracevara, \conj{\tracevara}\right\rangle_{S_R}\big)
		\qquad\text{for all }\tracevara\in H^{-1/2}(S_R)
	\end{align*}
	holds for its inverse $\NtD=\DtN^{-1}$ as well.
	This and the previous identities verify
	\begin{equation}\label{eqn:Re_b_bound}
		\begin{aligned}
			\Re b( (\mathbf{p},\nomathfrak{u}), (\opA\nabla {w}, {w}))
			&= \|u\|_{L^2(B_R)}^2,\\
			\Re b( (\mathbf{p},\nomathfrak{u}), (-\mathbf{p}, u))
			 &\geq \|\opA^{-1/2}\mathbf{p}\|_{L^2(B_R)}^2 
				 +\Re(s^2)\|p^{1/2}u\|_{L^2(B_R)}^2,\\
			-\Re b((\mathbf{p},\nomathfrak{u}), (0, s^2 u+p^{-1}\Div \mathbf{p}))
			 &\geq \|p^{-1/2}\Div\mathbf{p}\|_{L^2(B_R)}^2 
				 - |s|^4 \|p^{1/2}u\|_{L^2(B_R)}^2.
		\end{aligned}
	\end{equation}
	Define $\mathbf{q}\in H(B_R,\Div)$ and $v\in L^2(B_R)$ with 
	$c(s)\coloneqq\max\{0,p_{\max}-a_{\max}\Re(s^2)/|s|^2\}$ by
	\begin{align}\label{eqn:m_v_infsup}
		\begin{pmatrix}
			\mathbf{q}\\v
		\end{pmatrix}
		=a_{\max}\begin{pmatrix}
			-\mathbf{p}\\u
		\end{pmatrix}
		- p_{\max}|s|^{-2}\begin{pmatrix}
			0\\s^2 u +p^{-1}\Div\mathbf{p}
		\end{pmatrix}
		+(1+c(s)p_{\max})|s|^2\begin{pmatrix}
			\opA\nabla {w}\\{w}
		\end{pmatrix}.
	\end{align}
	The combination~\eqref{eqn:Re_b_bound}--\eqref{eqn:m_v_infsup} results in
	\begin{align}\nonumber
		\Re b&( (\mathbf{p},\nomathfrak{u}), (\mathbf{q}, v))
		\geq a_{\max}\|\opA^{-1/2}\mathbf{p}\|_{L^2(B_R)}^2
		+p_{\max}|s|^{-2}\|p^{-1/2}\Div\mathbf{p}\|_{L^2(B_R)}^2\\\nonumber
		 &\quad\;+ (1+ c(s)p_{\max})|s|^2\|u\|_{L^2(B_R)}^2
		 + (a_{\max}\Re(s^2)/|s|^2-p_{\max})|s|^2\|p^{1/2} u\|_{L^2(B_R)}^2\\
		  &\geq |s|^{-2}\|\Div\mathbf{p}\|_{L^2(B_R)}^2
			+ \|\mathbf{p}\|_{L^2(B_R)}^2
		 + |s|^2\|u\|_{L^2(B_R)}^2 
		 = \|(\mathbf{p},u)\|_{M,s}^2.\label{eqn:b_infsup}
	\end{align}
	Triangle inequalities for~\eqref{eqn:m_v_infsup},
	$|s|\|{w}\|_{L^2(B_R)}\leq \|{w}\|_{V(B_R,\opA,s)}$ by~\eqref{eqn:V_norm_def},
	and~\eqref{eqn:b_infsup_U_bound} reveal
	\begin{align*}
		\|\mathbf{q}\|_{H(B_R,\Div),s}
		&\leq a_{\max}\|\mathbf{p}\|_{H(B_R,\Div),s} +
		(1+c(s)p_{\max})|s|\tCN(s)\|u\|_{L^2(B_R)},\\
		\|v\|_{L^2(B_R)}
		&\leq \frac{p_{\max}}{p_{\min}}|s|^{-2}\|\Div \mathbf{p}\|_{L^2(B_R)}
		+\left(a_{\max}\!+\!p_{\max}\!+\!(1\!+\!c(s)p_{\max})\tCN(s)\right)\|u\|_{L^2(B_R)}.
	\end{align*}
	Since $c(s)\leq a_{\max}+p_{\max}$ by definition, the previous estimates and~\eqref{eqn:M_norm_def} establish
	\begin{align}\label{eqn:inf_sup_norm_estimate}
		\|(\mathbf{q},v)\|_{M,s}\leq 
		\Const{b}\big(1+\tCN(s)\big)\|(\mathbf{p},u)\|_{M,s}
	\end{align}
	for a constant %
	$\Const{b}>0$ that exclusively depends on 
	$a_{\max},p_{\max}$, and $p_{\min}$.
	This and~\eqref{eqn:b_infsup} provide the inf-sup condition
	\begin{align*}
		\inf_{0\ne(\mathbf{p},u)\in M}\sup_{0\ne(\mathbf{q},v)\in M}
			\frac{\Re b\left((\mathbf{p},u),(\mathbf{q},v)\right)}
			{\|(\mathbf{p},u)\|_{M,s}\|(\mathbf{q},v)\|_{M,s}}
			\geq \left( \Const{b}\big(1+\tCN(s)\big)\right)^{-1}>0.
	\end{align*}
	Analogical arguments with $(\mathbf{q},v)\in M$
	from~\eqref{eqn:m_v_infsup} with $s$ 
	replaced by $\conj s$
	reveal the inf-sup condition for the adjoint problem.
	Hence~\eqref{eqn:D_mixed_def} is well posed and admits a unique solution $(\mathbf{p},u)\in M$.
	
	\medskip
	\noindent\emph{Step 2} (characterisation):
	To verify the norm estimate \emph{(i)}, 
	employ~\eqref{eqn:D_mixed_def} and~\eqref{eqn:b_infsup} for
		\begin{align*}
			\|(\mathbf{p},u)\|_{M,s}^2
			\leq \Re b( (\mathbf{p},u),(\mathbf{q},v)) 
			= \Re \left\langle g,\traceNu[\IGebiet]\conj{\mathbf{q}}\right\rangle _{\Interface}.
		\end{align*}
	Observe 
	$\Const{tr,N}^{-1}\|\traceNu[\IGebiet]\mathbf{q}\|_{H^{-1/2}(\Interface),s}
	\leq \|\mathbf{q}\|_{H(\IGebietInt,\Div),s}\leq
	\|(\mathbf{q},v)\|_{M,s}$ from~\eqref{eqn:Ctr_N} 
	for the continuous normal trace
	$\traceNu[\IGebiet]\mathbf{q}=\traceNu[\IGebietInt]\mathbf{q}$ 
	of $\mathbf{q}\in H(B_R,\Div)$ and~\eqref{eqn:M_norm_def}.
	Hence \eqref{eqn:inf_sup_norm_estimate} and \Cref{lem:tCN_bound} result with 
	$\Const{DL}\coloneqq \Const{tr,N}\Const{b}\big(1+\sqrt{2+(2p_{\max}^2 + 1)}\big)$~in
		\begin{align*}
			\|(\mathbf{p},u)\|_{M,s}
			\leq \Const{b}(1+\tCN(s))\|g\|_{H^{1/2}(\Interface),s}
			\leq \Const{DL}(1+\CN(s))\|g\|_{H^{1/2}(\Interface),s}.
		\end{align*}
	This is \emph{(i)} and it remains to prove \emph{(ii)--(iv)}.
	Consider any $\tracevara\in H^{1/2}(\Interface)$ and the unique solution 
	$(\mathbf{p},u)\in M$ to \eqref{eqn:D_mixed_def}.
	The mixed problem
	\eqref{eqn:D_mixed_def} for 
	$\mathbf{q}\coloneqq \opA \nabla v\in H(B_R,\Div)$
	reveal with
	an integration by parts (as in~\eqref{eqn:jv_ibp} for $v$ instead of ${w}$) that
	\begin{align*}
		\left\langle \tracevara,
		\traceN[\IGebiet]\conj{v}\right\rangle_{\Interface}
		=b\left((\mathbf{p},u),(\mathbf{q},v)\right)
		= \left\langle \nomathfrak{u}, \opL\conj{v}\right\rangle_{B_R}
		\qquad\text{for all } v\in V(B_R,\opA,\conj s).
	\end{align*}
	This proves \emph{(ii)}.
	The first equation of~\eqref{eqn:D_mixed_def} and
	$\traceNu[\IGebiet]\mathbf{q}=0=\traceNu[B_R]\mathbf{q}$ for all
	$\mathbf{q}\in C^\infty_0(B_R\setminus\Interface;\R^\dim)$ 
	implies that $\opA^{-1}\mathbf{p} = \nabla \nomathfrak{u}$ is the weak gradient of
	$\nomathfrak{u}$ in $L^2(B_R\setminus\Interface)$.
	In other words, 
	$u\in H^1(B_R\setminus\Interface,\opA)$ holds with
	$\mathbf{p}|_{B_R\setminus\Interface}=\opA\nabla u|_{B_R\setminus\Interface}
	\in H(B_R\setminus\Interface,\Div)$.

	It remains to prove $\jumpNext u=0$ for \emph{(iii)} and 
	the jump relations \emph{(iv)}.
	Let $\mathbf{q}\in H(B_R,\Div)$ be arbitrary and observe
	$\left\langle\NtD u,\conj{\traceNu[B_R]\mathbf{q}}\right\rangle_{S_R}
	=\big\langle \traceN[B_R]u,
	\conj{\NtD[\conj s] \traceNu[B_R]\mathbf{q}}\big\rangle_{S_R}$
	from the corresponding identity for the inverse $\DtN$ as in the proof of 
	\Cref{lem:Newton_potential_L2}.
	Since $\Interface$ has measure zero, 
	this and an integration by parts over 
	$B_R\setminus\Interface$ 
	with the first equation of~\eqref{eqn:D_mixed_def}
	and $\mathbf{p}|_{B_R\setminus\Interface}
	=\opA\nabla u|_{B_R\setminus\Interface}\in H(B_R\setminus\Interface,\Div)$
	from \emph{(iii)}
	verify
	\begin{align*}
		\left\langle \tracevara, \conj{\traceNu[\IGebiet]\mathbf{q}}\right\rangle_{\Interface}
		&=-\int_{B_R\setminus\Interface}\left(\nabla u\cdot\conj{\mathbf{q}} + u\,
			\Div\conj{\mathbf{q}}\right) \d x
		+\big\langle \traceN u, 
		\conj{\NtD[\conj s]\traceNu[B_R]\mathbf{q}}\big\rangle_{S_R}\\
		&=-\big\langle \jumpD[\IGebiet]{u},
		\conj{\traceNu[\IGebiet]\mathbf{q}}\big\rangle_{\Interface}
		+\big\langle \NtD\jumpNext u, 
		\conj{\traceNu[B_R]\mathbf{q}}\big\rangle_{S_R}
	\end{align*}
	with $-\NtD\jumpNext u = (\traceD -\NtD) u$ in the last step.
	Since the boundaries $\Interface$ and $S_R$ are separated 
	($\mathrm{dist}(\Interface,S_R)>0$ by $\Interface\subset B_R$),
	the normal components of functions in $H(B_R,\Div)$ are 
	independent and surjective onto $H^{-1/2}(\Interface)\times H^{-1/2}(S_R)$.
	Hence the previous identity and the injectivity of $\NtD$ verify
	$\jumpNext u =0$, implying \emph{(iii)}
	and 
	$\jumpD[\IGebiet] u= -\tracevara$.
	This and $\jumpN[\IGebiet]u=0$ from
	the continuity of the normal component
	$\traceN[\IGebiet] u=\traceNu[\IGebiet]\mathbf{p} = -\traceNext[\IGebiet] u$
	of $\mathbf{p}\in H(B_R,\Div)$ across $\Interface$ by \emph{(iii)} reveal
	\emph{(iv)} and conclude the proof.
\end{proof}

\begin{theorem}[double layer potential]\label{lem:D_properties}
	The double layer potential $\opD$ 
	from~\eqref{eqn:D_def}
	maps $H^{1/2}(\Interface)$ boundedly into 
	$L^2(B_R)\cap V(B_R\setminus\Interface,\opA,s)$ 
	and is uniquely defined by~\eqref{eqn:D_ultraweak_def}.
	Any $g\in H^{1/2}(\Interface)$ and $u\coloneqq\opD g$ satisfy
	with the constant $\Const{DL}>0$ from \Cref{lem:well-posedness of mixed formulation}	that
	\begin{enumerate}[label=(\roman*)]
		\item $\displaystyle
		\|u\|_{V(B_R\setminus\Interface,\opA, s)}
		\leq \Const{DL}\left(1+\CN(s)\right)\|g\|_{H^{1/2}(\Interface),s}$,
		\item $-\Div(\opA \nabla u) + s^2 p\, u = 0$ in $B_R\setminus\Interface$, and
		\item $
		\jumpD[\IGebiet]{u} = -\tracevara
		\quad\text{and}\quad
		\jumpN[\IGebiet]{u} = 0
$.
	\end{enumerate}
\end{theorem}
\begin{proof}[Proof of \Cref{lem:D_properties}]
	The boundedness of $\opNdual:L^2(B_R)\to V(B_R,\opA,\conj s)$ 
	by \Cref{lem:Newton_potential_L2}
	implies the boundedness of its adjoint
	$\opN=\opNddual:(V(B_R,\opA,\conj s))'\to L^2(B_R)$ defined by%
	~\eqref{eqn:Nddual_def}.
	Hence $\opD=\opN\circ \traceN[\IGebiet]'$ is a bounded operator.
	\Cref{lem:Newton_potential_L2}
	provides $\opL\conj v =\conj{\opLdual v}$ for all $v\in V(B_R,\opA,\conj s)$
	and the definition~\eqref{eqn:D_def}
	of $\opD$ results~in
	\begin{align*}
		\left\langle \opD\tracevara, \opL \conj{v}\right\rangle_{B_R}
		=\big\langle \tracevara, \traceN[\IGebiet]\conj{\opNdual\opLdual v}\big\rangle_{\Interface}
		\quad\text{for all }\tracevara\in H^{1/2}(\Interface), v\in V(B_R,\opA,\conj s).
	\end{align*}
	This and $\opNdual\opLdual=\id$ verify~\eqref{eqn:D_ultraweak_def}.
	Since $\opL:V(B_R,\opA,s)\to L^2(B_R)$ is surjective, \eqref{eqn:D_ultraweak_def}
	uniquely defines $\opD\tracevara$.
	Consider any $\tracevara\in H^{1/2}(\Interface)$ and set $u\coloneq\opD\tracevara$.
	For any $\oldvarphi\in C^\infty_0(B_R\setminus\Interface)$, 
	\eqref{eqn:D_ultraweak_def} with $\traceN[\IGebiet] \conj \oldvarphi=0$ and
	the characterisation of $\opL \conj{\oldvarphi}$ by~\eqref{eqn:L_L2_repr} show
	\begin{align*}
		0=\left\langle u, \opL \conj{\oldvarphi}\right\rangle_{B_R} 
		= \int_{B_R}u\,(-\Div(\opA\nabla \conj{\oldvarphi})) \d x +
		\int_{B_R}s^2p\, u\, \conj\oldvarphi\d x.
	\end{align*}
	The definition of weak derivatives reveals
	$-\Div(\opA\nabla u) = -s^2p\, u\in L^2(B_R\setminus\Interface)$, 
	implying~$(ii)$.
	The characterisation of the unique solution $(\mathbf{p},u)\in M$ 
	to~\eqref{eqn:D_mixed_def} in
	\Cref{lem:well-posedness of mixed formulation}.ii--iii 
	establishes $\opD \tracevara=u\in V(B_R\setminus\Interface,\opA,s)$
	and the jump relations \emph{(iii)}.
	Observe
	\begin{align*}
		\|\opD g\|_{V(B_R\setminus\Interface,\opA, s)}\leq\|(\mathbf{p},u)\|_{M,s}
		\leq \left(1+\CN(s)\right)\Const{DL}\|g\|_{H^{1/2}(\Interface),s}
	\end{align*}
	from the definition of the
	involved norms with $\opA\nabla u=\mathbf{p}\in L^2(B_R\setminus\Interface)$ by \Cref{lem:well-posedness of mixed
	formulation}.i and \Cref{lem:well-posedness of mixed formulation}.iv.
	This concludes the proof.
\end{proof}

\subsection{The Calder\'on operator}%
\label{sub:The Calder'on operator}
The transmission problem~\eqref{eqn:IP} and its reformulation~\eqref{eqn:TIP}
on the truncated domain $B_R\setminus\Interface$ prescribe Dirichlet and Neumann jumps 
across the interface $\Interface$ and its solutions are characterised by 
the single and double layer operators 
from \Cref{sub:Single_layer_potential,sub:Double_layer_potential}.
\begin{lemma}[representation formula]\label{lem:integral_formulation}
	Given any
	$\gD\in H^{1/2}(\Interface)$ and $\gN\in H^{-1/2}(\Interface)$,
	the unique solution $u\in V(B_R\setminus\Interface,\opA,s)$ to~\eqref{eqn:TIP} 
	reads
	\begin{align}\label{eqn:interface_integral_formulation}
		u = \opS \gN - \opD \gD.
	\end{align}
	In particular, any $v\in V(B_R\setminus\Interface,\opA,s)$ with $-\Div(\opA\nabla v) + s^2p\, v=0$ satisfies
	\emph{Green's representation formula}
	\begin{align}\label{eqn:Green_representation}
		v = \opS\jumpN[\IGebiet] v - \opD\jumpD[\IGebiet] v.
	\end{align}
\end{lemma}
\begin{proof}
	Any solution $u$ to~\eqref{eqn:TIP} for $\gD=0=\gN$ satisfies 
	$u\in V(B_R,\opA,s)$ and solves~\eqref{eqn:THP} (for $\Omega_R=B_R$).
	The uniquess follows from \Cref{ass:A_p_new}.
	The jump relations from the characterisation of
	$\opS$ and $\opD$ in \Cref{lem:S_property,lem:D_properties} below
	verify~\eqref{eqn:interface_integral_formulation}.
\end{proof}

Green's representation formula~\eqref{eqn:Green_representation} enables
a reformulation of the transmission problem~\eqref{eqn:TIP} as boundary (integral) equations for the Cauchy traces of
solutions on the interface $\Interface$.
Since the jumps are prescribed by the transmission problem, 
the Cauchy traces are uniquely defined by the 
averages $\meanD[\IGebiet]\bullet$ and $\meanN[\IGebiet]\bullet$ 
from~\eqref{eqn:jump_mean_def}.
The maps
\begin{align}\label{eqn:Bop_V}
	\opV\phantom{'} &: H^{-1/2}(\Interface)\to H^{1/2}(\Interface)
		 &&\text{with}\qquad
		 \opV \gN\coloneqq\meanD[\IGebiet]{\opS \gN},\\\label{eqn:Bop_K}
	\opK\phantom{'} &: H^{1/2}(\Interface)\to H^{1/2}(\Interface)
		 &&\text{with}\qquad
		 \opK \gD\coloneqq \meanD[\IGebiet]{\opD \gD},\\\label{eqn:Bop_Kdual}
	\opKdual &: H^{-1/2}(\Interface)\to H^{-1/2}(\Interface)
			 &&\text{with}\qquad
			 \opKdual \gN\coloneqq \meanN[\IGebiet]{\opS \gN},\\\label{eqn:Bop_W}
	\opW &: H^{1/2}(\Interface)\to H^{-1/2}(\Interface)
		 &&\text{with}\qquad
		 \opW \gD\coloneqq \meanN[\IGebiet]{\opD \gD}
\end{align}
for any $\gD\in H^{1/2}(\Interface)$ and $\gN\in H^{-1/2}(\Interface)$ 
are called \emph{single layer}, \emph{double
layer}, \emph{dual double layer}, and \emph{hypersingular boundary integral operators}, respectively.
Recall $\CN(s)$ from~\eqref{eqn:CN_def}, $\Const{tr,N}$ from \Cref{lem:trace_inequality}
as well as $\Const{SL}$ and $\Const{DL}$
from \Cref{lem:S_property} and \Cref{lem:well-posedness of mixed formulation}.
\begin{lemma}[boundedness]\label{lem:boundedness}
	The boundary operators~\eqref{eqn:Bop_V}--\eqref{eqn:Bop_W}
	are bounded with
	\begin{align*}
		\|\opV\gN\|_{H^{1/2}(\Interface),s}
		&\leq \CN(s)\,\|\gN\|_{H^{-1/2}(\Interface),s},\\
		\|\opK\gD\|_{H^{1/2}(\Interface),s}
		&\leq 
		\left(\tfrac12+\Const{DL}(1+\CN(s))\right)\,\|\gD\|_{H^{1/2}(\Interface),s},
		\\
		\|\opKdual\gN\|_{H^{-1/2}(\Interface),s}
		&\leq \left(\tfrac12+\Const{tr,N}\Const{SL}\CN(s)\right)
		\|\gN\|_{H^{-1/2}(\Interface),s},\\
		\|\opW\gD\|_{H^{-1/2}(\Interface),s}
		&\leq \Const{tr,N}\Const{DL}(1+\CN(s))\,\|\gD\|_{H^{1/2}(\Interface),s}
	\end{align*}
	for all $\gD\in H^{1/2}(\Interface)$ and $\gN\in H^{-1/2}(\Interface)$.
\end{lemma}
\begin{proof}
	By the jump relations in \Cref{lem:S_property,lem:D_properties},
	the Dirichlet (resp.~Neumann) traces of $\opS$ (resp.~$\opD$) are single-valued
	such that $\opV=\traceD[\IGebietR] \opS$ and $\opW=\traceN[\IGebietInt]\opD$.
	Hence \Cref{lem:trace_inequality} and \Cref{lem:S_property}.i reveal 
	for any $\gN\in H^{-1/2}(\Gamma)$ that
	\begin{align*}
		\|\opV\gN\|_{H^{1/2}(\Interface),s}
		\leq \|\opS\gN\|_{H^{1}(\IGebietR),s}
		\leq \CN(s)\|\gN\|_{H^{-1/2}(\Interface),s}.
	\end{align*}
	Similarly, \Cref{lem:trace_inequality} 
	(for $\mathbf{p}=\opA\nabla \opD\gD$) and \Cref{lem:D_properties}.i lead for
	any $\gD\in H^{1/2}(\Interface)$~to
	\begin{align*}
		\Const{tr,N}^{{-1}}\|\opW\gD\|_{H^{-1/2}(\Interface),s}
		&\leq \|\opD\gD\|_{H(\IGebietInt,\Div),s}
		\leq \Const{DL}(1+\CN(s))\|\gD\|_{H^{1/2}(\Interface),s}.
	\end{align*}
	This proves the claimed bounds for $\opV$ and $\opW$.
	The combination of
	\Cref{lem:trace_inequality} (with $\mathbf{p}=\opA\nabla\opS \gN$)
	with Theorem
	\ref{lem:S_property}.i and \ref{lem:D_properties}.i verify as before that
	\begin{align*}
		\|\traceD[\IGebiet]{\opD\gD}\|_{H^{1/2}(\Gamma),s}
		&\leq\|{\opD\gD}\|_{H^{1}(\IGebietR),s}
		\leq \Const{DL}(1+\CN(s)) \|\gD\|_{H^{1/2}(\Gamma),s},\\
		\Const{tr,N}^{-1}\|\traceN[\IGebietInt]{\opS\gN}\|_{H^{-1/2}(\Gamma),s}
		&\leq\|\opA\nabla{\opS\gN}\|_{H(\IGebietInt,\Div),s}
		\leq \Const{SL}\CN(s) \|\gN\|_{H^{-1/2}(\Gamma),s}
	\end{align*}
	for all $\gD\in H^{1/2}(\Gamma)$ and $\gN\in H^{-1/2}(\Gamma)$.
	This and triangle inequalities with
	\begin{align*}
		\left|\traceD[\IGebietR]{\opD\gD} - \opK\gD \right|
		&%
		= \left|\tfrac{1}{2}\jumpD[\IGebietR]{\opD\gD}\right|
		= \left|\tfrac{1}{2}\gD\right|,\\
		\left|\traceN[\IGebietInt]{\opS\gN} - \opKdual\gN \right|
		&%
		= \left|\tfrac{1}{2}\jumpN[\IGebietInt]{\opS\gN}\right|
		= \left|\tfrac{1}{2}\gN\right|,
	\end{align*}
	using~\eqref{eqn:jump_mean_def} and 
	the jump relations (with respect to $\IGebietR$) 
	from~\Cref{lem:S_property,lem:D_properties},
	provide the remaining bounds and conclude the proof.
\end{proof}

The \emph{Calder\'on operator} $\opC$ 
on the Cauchy trace space $\mathbf{X}(\Interface)\coloneqq H^{1/2}(\Interface)\times H^{-1/2}(\Interface)$ reads
\begin{align}\label{eqn:Calderon_def}
	\opC\coloneqq
	\begin{pmatrix}
		-\opK&\opV\\-\opW&\opKdual
	\end{pmatrix}:\mathbf{X}(\Interface)\to\mathbf{X}(\Interface).
\end{align}
By
the Green representation formula~\eqref{eqn:Green_representation},
any $\boldsymbol{g}= (\gD,\gN)\in \mathbf{X}(\Interface)$ defines
a (unique) solution
$u\in V(B_R\setminus\Interface,\opA,s)$ to~\eqref{eqn:TIP} with
\begin{align}\label{eqn:g_ident}
	\boldsymbol{g} = \big(\jumpD[\IGebiet] u,\jumpN[\IGebiet]u\big)
	\qquad\text{and}\qquad
	\opC \boldsymbol{g} 
	= \big(
	\meanD[\IGebiet]{u},\meanN[\IGebiet]{u}
	\big).
\end{align}
This and the definition of the jumps and averages in~\eqref{eqn:jump_mean_def} imply the \emph{Calder\'on identity}.
\begin{lemma}[Calder\'on identity]\label{lem:C_equiv}
	Any $\boldsymbol{g} =(\gD,\gN)\in \mathbf{X}(\Interface)$
	and the unique solution $u\in V(B_R\setminus\Interface,\opA,s)$ to~\eqref{eqn:IP} 
	(with $u = \opS\gN-\opD\gD$ by \Cref{lem:integral_formulation}) satisfy
	\begin{align*}
		\big(\opC+\tfrac12\big)\boldsymbol{g} 
			=\big(\traceD[\IGebiet]u,\traceN[\IGebiet] u\big)
		\qquad\text{and}\qquad
		\big(\opC-\tfrac12\big)\boldsymbol{g} 
			=\big(\traceDext[\IGebiet]u,-\traceNext[\IGebiet] u\big).
	\end{align*}
\end{lemma}
\noindent In particular, \Cref{lem:C_equiv} (whose proof is omitted) reveals the equivalences%
\begin{align*}
	\opC \boldsymbol{g} = \tfrac{1}{2}\boldsymbol{g}
	\;\Leftrightarrow\;
	\boldsymbol{g} = \big(\traceD[\IGebiet]u,\traceN[\IGebiet] u\big),
	\qquad
	\opC \boldsymbol{g} = -\tfrac{1}{2}\boldsymbol{g}
	\;\Leftrightarrow\;
	\boldsymbol{g} = \big(\traceDext[\IGebiet]u,-\traceNext[\IGebiet] u\big).
\end{align*}
The Cauchy trace space $\mathbf{X}(\Interface)$ equipped with the usual product norm 
is self-dual with the dual pairing given for any
	$\boldsymbol{g}=(g_D,g_N),\boldsymbol{h}=(h_D,h_N)\in \mathbf{X}(\Interface)$ by
\begin{align}\label{eqn:X_dual}
	\left\langle \boldsymbol{g},\boldsymbol{h}\right\rangle_{\mathbf{X}(\Interface)}
	\coloneqq
	\left\langle g_D,h_N\right\rangle _{\Interface}
	+
	\left\langle h_D, g_N\right\rangle _{\Interface}.
\end{align}
The induced norm reads
\begin{align}\label{eqn:X_norm}
	\|\boldsymbol{g}\|_{\mathbf{X}(\Interface),s}
	\coloneqq\sqrt{\NormHD[\Interface]{s}{\gD}^2+\NormHN[\Interface]{s}{\gN}^2}
	\qquad\text{for all }\boldsymbol{g} =(\gD,\gN)\in\mathbf{X}(\Interface).
\end{align}
Recall $\slow$ and $\sLow=\slow^{-1}$ from~\eqref{eqn:s_low_def}.
\begin{lemma}[Calder\'on operator]\label{lem:C_garding}
	There exists a compact operator 
	$\mathsf T(s):\mathbf{X}(\Interface)\to \mathbf{X}(\Interface)$ such that 
	$\opC:\mathbf{X}(\Interface)\to\mathbf{X}(\Interface)$ from~\eqref{eqn:Calderon_def} 
	and any $\boldsymbol{\tracevarb} \in \mathbf{X}(\Interface)$ satisfy
	\begin{align*}
		\|\opC \boldsymbol{g} \|_{\mathbf{X}(\Interface),s}
		&\leq \Const{G}(1+\CN(s))\,
		\|\boldsymbol{g}\|_{\mathbf{X}(\Interface),s},\\
		\Re\left\langle (\opC+\mathsf T(s)) \boldsymbol{\tracevarb} ,
		\conj{\boldsymbol{\tracevarb}} \right\rangle_{\mathbf{X}(\Interface)}
		&\geq c_{\mathrm{G}}^2 \,\slow^2\,
		\|\boldsymbol{\tracevarb} \|_{\mathbf{X}(\Interface),s}^2
	\end{align*}
	The constants $c_{\mathrm{G}},\Const{G}>0$ are independent of $s$ and 
	exclusively depend on $\Interface,R,\opA$, and $p$.
\end{lemma}
\begin{proof}
	The boundedness of $\opC$ follows from~\eqref{eqn:Calderon_def},~\eqref{eqn:X_norm},
	and \Cref{lem:boundedness}.
	Consider
	any $\boldsymbol{\tracevarb} =(\tracevarb_D,\tracevarb_N)\in \mathbf{X}(\Interface)$ and
	set $u\coloneqq \opS\tracevarb_{\mathrm{N}} - \opD\tracevarb_{\mathrm{D}}\in V(B_R\setminus\Interface,\opA,s)$.
	To prove the coercivity of $\opC+\mathsf T(s)$ for some compact operator 
	$\mathsf T(s)$, we first establish
	\begin{align}\label{eqn:RC_xi_identity}
		\Re\left\langle \opC \boldsymbol{\tracevarb} ,\conj{\boldsymbol{\tracevarb} }\right\rangle_{\mathbf{X}(\Interface)} 
		= \|\opA^{1/2}\nabla u\|_{L^2(B_R\setminus\Interface)}^2 + \Re(s^2)\|p^{1/2}u\|_{L^2(B_R)}^2 
		- \Re\left\langle \DtN u,\conj u\right\rangle_{S_R}.
	\end{align}
	Elementary algebra,~\eqref{eqn:g_ident}, and the product rule 
	$\jump{{A} {B}}=\jump{{A}}\mean{{B}}+\mean{{A}}\jump{{B}}$
	for jumps~show%
	\footnote{The product rule for jumps and 
		$\mean{\conj {B}}\jump{{A}}= \conj{\mean{{B}}\jump{\conj {A}} }$
		imply
		$\mean{{A}}\jump{\conj {B}} + \mean{{B}}\jump{\conj {A}}
		= [{A}\conj {B}] + 2i\Im\mean{{B}}\jump{\conj {A}}$.
	}
	\begin{align*}
		&\left\langle \opC\boldsymbol{\tracevarb}  ,
		\conj{\boldsymbol{\tracevarb} }\right\rangle_{\mathbf{X}(\Interface)}
		=\big\langle\meanD[\IGebiet]{u},\jumpN[\IGebiet]{\conj{u}}\big\rangle_{\Interface}
		+ \big\langle \meanN[\IGebiet] {u},\jumpD[\IGebiet]{\conj u}\big\rangle_{\Interface}
		\\
		&\qquad=\int_{\Interface}\left(\traceD[\IGebiet]{u}\traceN[\IGebiet]{\conj{u}} + 
			\traceDext[\IGebiet]{u}\traceNext[\IGebiet]{\conj{u}}\right) \d s
		+ 2i\Im \big\langle \meanN[\IGebiet]{u},\jumpD[\IGebiet]{\conj u}
		\big\rangle_{\Interface}.
	\end{align*}
	A piecewise integration by parts on $\IGebiet$ and $B_R\setminus\overline{\IGebiet}$ 
	results for the
	real part in
	\begin{align*}
		&\Re\left\langle \opC\boldsymbol{\tracevarb}  ,
		\conj{\boldsymbol{\tracevarb} }\right\rangle_{\mathbf{X}(\Interface)}
		=\Re\int_{B_R\setminus\Interface}\left(\opA\nabla u\cdot\nabla\conj u + \Div(\opA\nabla u)\,\conj u\right)\d x
		-\Re\left\langle \traceN[B_R] u, \conj{u}\right\rangle_{S_R}.
	\end{align*}
	Recall $\Div(\opA\nabla u) = s^2p u$ from \Cref{lem:S_property,lem:D_properties} so that
	$\jumpNext u = 0$ from $u\in V(B_R\setminus\Interface,\opA,s)$ proves the 
	identity~\eqref{eqn:RC_xi_identity}.
	Moreover,~\ref{ass:C2} and
	$\Div(\opA\nabla u)=s^2p u$ lead to
	\begin{align*}
		\|u\|_{H^1(B_R\setminus\Interface),s}^2
		&\leq \max\{a_{\min}^{-1},p_{\min}^{-1}\}
		\left(\|\opA^{1/2} \nabla u\|_{L^2(B_R\setminus\Interface)}^2 + 
		|s|^2\|p^{1/2}u\|_{L^2(B_R)}^2\right),\\
		\|\opA\nabla u\|_{H(B_R\setminus\Interface,\Div),s}^2
		&\leq \max\{a_{\max},p_{\max}\}
		\left(\|\opA^{1/2} \nabla u\|_{L^2(B_R\setminus\Interface)}^2 + 
		|s|^2\|p^{1/2}u\|_{L^2(B_R)}^2\right).
	\end{align*}
	The $s$-explicit trace inequality of \Cref{lem:trace_inequality} and 
	a Cauchy inequality provide %
	\begin{align*}%
		\|\traceD[\IGebiet]{u}\|_{H^{1/2}(\Interface),s} 
			+ \|\traceDext[\IGebiet]{u}\|_{H^{1/2}(\Interface),s}
		&\leq (1+\Const{tr,D}^2)^{1/2} \sLow\|{u}\|_{H^1(B_R\setminus\Interface),s},\\
		\|\traceN[\IGebiet]{u}\|_{H^{-1/2}(\Interface),s} 
		+ \|\traceNext[\IGebiet]{u}\|_{H^{-1/2}(\Interface),s}
			&\leq 2\Const{tr,N}\sLow^2
			\|\opA\nabla {u}\|_{H(B_R\setminus\Interface,\Div),s}.
	\end{align*}
	Hence triangle inequalities with 
	$\boldsymbol{g} =(\jumpD[\IGebiet]u,\jumpN[\IGebiet]u)$
	from~\eqref{eqn:g_ident} %
	imply
	\begin{align*}
		\|\boldsymbol{\tracevarb} \|_{\mathbf{X}(\Interface)}^2
		\leq \Const{tr}\,\sLow^2 %
			\left(\|\opA^{1/2}\nabla u\|_{L^2(B_R\setminus\Interface)}^2 + 
			|s|^2\|p^{1/2}u\|_{L^2(B_R)}^2\right)
	\end{align*}
	for a constant $\Const{tr}>0$ that exclusively depends on 
	$\Const{tr,D},\Const{tr,N},a_{\max},a_{\min},p_{\max}$, and $p_{\min}$.
	Consequently,~\eqref{eqn:RC_xi_identity} with $\Re (s^2)\leq |s|^2$
	and $-\Re\left\langle \DtN u,\conj u\right\rangle_{S_R}\geq 0$ 
	from~\eqref{eqn:DtN_non_pos}
	reveal
	\begin{align}\label{eqn:C_garding_final}
		\Const{tr}^{-1}\,
		\slow^2\|\boldsymbol{\tracevarb} \|_{\mathbf{X}(\Interface)}^2\leq 
		\Re\left\langle\opC\boldsymbol{\tracevarb},
			\conj{\boldsymbol{\tracevarb}}\right\rangle_{\mathbf{X}(\Interface)}
		+ 2|s|^2\|p^{1/2}u\|_{L^2(B_R)}^2.
	\end{align}
	Since $\opS$ and $\opD$ map boundedly into 
	$V(B_R\setminus\Interface,\opA,s)
	\subset H^1(B_R\setminus\Interface)\hookrightarrow L^2(B_R)$
	and the latter embedding is compact \cite[Thm.~4.11]{EG:MeasureTheoryFine2015}, 
	the map $\mathsf T(s):\mathbf{X}(\Interface)\to \mathbf{X}(\Interface)$
	given~by
	\begin{align*}
		\left\langle \mathsf T(s) \boldsymbol{g}, 
		\conj{\boldsymbol{h}} \right\rangle_{\mathbf{X}(\Interface)}
			=2|s|^2\,\int_{B_R}p\,\left(\opS g_{N}-\opD g_D\right)\,
				\big( \conj{\opS h_{N}-\opD h_D}\big)\d x
	\end{align*}
	for any 
	$\boldsymbol{g} =(g_{D},g_{N}),\boldsymbol{h} =(h_{D},h_{N})\in \mathbf{X}(\Interface)$
	is a bounded compact operator.
	Since
	$\left\langle \mathsf T(s)\boldsymbol{\tracevarb},
	\conj{\boldsymbol{\tracevarb} }\right\rangle_{\mathbf{X}(\Interface)}
	=2|s|^2\|p^{1/2}u\|_{L^2(B_R)}^2$,%
	~\eqref{eqn:C_garding_final} concludes the proof with 
	$c_{\mathrm{G}}\coloneqq \Const{tr}^{-1/2}$.
\end{proof}
\newpage
\section{Stable integral formulation of transmission problems}%
\label{sec:Stable integral formulation for the transmission problem}
Green's formula and the boundary layer operators from~\Cref{sec:Potential operators for interface problems} 
enable a stable formulation of transmission problems as skeleton integral equations (SIE)
for the Cauchy~data.

\subsection{The acoustic transmission problem}%
\label{ssub:The acoustic transmission problem}
The computational Lipschitz domain $\Omega\subset\R^n$ of the Helmholtz transmission 
problem~\eqref{eqn:ATP_intro} is the complement of the
bounded \emph{acoustic obstacle} $\R^n\setminus\Omega$
and partitioned into $J\in \N$ pairwise disjoint Lipschitz sets 
$\Omega_1,\dots,\Omega_J \subset \Omega$ and the unbounded component
\begin{align*}
	\Omega_0=\Omega\setminus\bigcup_{j=0}^J\overline{\Omega_j}
\end{align*}
as displayed in \Cref{fig:acoustic_domain}. 
The boundary $%
\partial\Omega=\Gamma_D\cup\Gamma_N$ of the acoustic obstacle $\R^\dim\setminus{\Omega}$
splits disjointly into the relatively closed Dirichlet ($\Gamma_D$)
and Neumann ($\Gamma_N$) parts.
Let 
\begin{align*}
	\Sigma\coloneqq \Gamma_0\cup\Gamma_1\cup\dots\cup\Gamma_J
	\qquad\text{with}\quad \Gamma_j\coloneq\partial\Omega_j,j=0,\dots,J
\end{align*}
denote the full transmission interface.
The wavenumber $s\in\Cstar_{\geq0}$ and the coefficients $\opA\in L^\infty(\Omega;\mathbb{S}^\dim)$
and $p\in L^\infty(\Omega)$ satisfy~\ref{ass:C1}--\ref{ass:C3} and \Cref{ass:A_p}
for a sufficiently large ball $B_R$ 
that also contains the transmission interface $\Sigma\subset B_R$, i.e.,
$\overline\Omega_1,\dots,\overline\Omega_J\subset B_R$.
Associated to the transmission interface are the multi-trace space%
~\cite{CHJ:MultitraceBoundaryIntegral2013,CHJP:NovelMultitraceBoundary2015}
\begin{align*}
	\Xmult(\hSigma)\coloneqq \prod_{j=0}^J \mathbf{X}(\Gamma_j)
	\qquad\text{with}\qquad
	\mathbf{X}(\Gamma_j)\coloneqq H^{1/2}(\Gamma_j)\times H^{-1/2}(\Gamma_j), j=0,\dots,J
\end{align*}
and the single-trace spaces (with and without boundary conditions for $\Gamma_{\mathrm{D}}$ 
and $\Gamma_{\mathrm{N}}$)
\begin{align}\label{eqn:X_single}
	\mathbf{X}(\Sigma)
		&\coloneqq\left\{(\traceD[\Omega_j] v, \traceNu[\Omega_j]\mathbf{q})_{j=0}^J\ :\ 
		v\in H^1(\Omega), \mathbf{q}\in H(\Omega,\Div)\right\}\subset\mathbb{X}(\Sigma),\\
	\mathbf{X}_0(\Sigma)
		&\coloneqq\left\{
			(\traceD[\Omega_j] v, \traceNu[\Omega_j]\mathbf{q})_{j=0}^J\ :\ 
			v\in H^1_{\Gamma_{\mathrm{D}}}(\Omega), 
		\mathbf{q}\in H_{\Gamma_{\mathrm{N}}}(\Omega,\Div)\right\}\subset\mathbb{X}(\Sigma).
		\label{eqn:X_single_hom}
\end{align}
The self-dual pairing and the weighted norm on $\mathbb{X}(\Sigma)$ inherited
from~\eqref{eqn:X_dual}--\eqref{eqn:X_dual} read
\begin{align}\label{eqn:self_dual_XX_def}
	\left\langle \bbg[],\bbh[]\right\rangle _{\mathbb{X}(\Sigma)}
	&\coloneqq \sum^{J}_{j=0} 
	\left\langle \bbg[j],\bbh[j]\right\rangle _{\mathbf{X}(\Gamma_j)}
	&&\text{for all }
	\bbg[]=(\bbg[j])_{j=0}^J,\bbh[]=(\bbh[j])_{j=0}^J\in\mathbb{X}(\Sigma),\\
\label{eqn:X_Sigma_norm}
	\|\boldsymbol{g}\|_{\mathbb{X}(\Sigma),s}
	&\coloneqq
	\sqrt{\sum^{J}_{j=0} 
	\|\boldsymbol{g}_j\|_{\mathbf{X}(\Interface_j),s}^2}
	&&\text{for all }\boldsymbol{g} =(\boldsymbol{g}_j)_{j=0}^J\in\mathbb{X}(\Sigma).
\end{align}%

The transmission condition~\eqref{eqn:ATP_intro_b}--\eqref{eqn:ATP_intro_e}
can be rewritten as $\traceC u - \bbg\in \mathbf{X}_0(\Sigma)$
with the Cauchy trace operator 
$\traceC:H^1_{\mathrm{loc}}(\Omega\setminus\Sigma,\opA)\to \Xmult(\hSigma)$ given by
\begin{align*}
\traceC v \coloneqq (\traceD[\Omega_j] v, \traceN[\Omega_j]v)_{j=0}^J\in\Xmult(\hSigma)
\qquad\text{for all }v\in H^1_{\mathrm{loc}}(\Omega\setminus\Sigma,\opA)
\end{align*}
and some $\bbg\in\mathbf{X}(\Sigma)$ that represents the Dirichlet and Neumann data
in~\eqref{eqn:ATP_intro_d}--\eqref{eqn:ATP_intro_e}.
To further allow inhomogeneities at the interior interfaces 
$\Gamma_j\cap \Gamma_k$ for $j,k=0,\dots,J$, 
we consider more general transmission data $\bbg\in\mathbb{X}(\Sigma)$.
The corresponding transmission problem~\eqref{eqn:ATP_intro} reads: 
Given $\bbg\in\mathbb{X}(\Sigma)$, find
a solution $u\in H^{1}_{\mathrm{loc}}(\Omega\setminus\Sigma)$ to
\begin{equation}\label{eqn:ATP}
	\begin{aligned}
		-\Div(\opA \nabla u) +s^2 p\, u = 0 
			\quad&\;\text{in }\Omega\setminus\Sigma,\\
			\traceC u-\bbg&\in \mathbf{X}_0(\Sigma),\\
	  \span \span u\text{ satisfies }~\eqref{radicond}.
	\end{aligned}
\end{equation}

To derive a boundary (integral) formalism of the transmission problem~\eqref{eqn:ATP},
we introduce an equivalent formulation in terms of Calder\'on operators 
(from \Cref{sub:The Calder'on operator})
for the (in general multi-valued) Cauchy traces 
$\traceC u\in\mathbb{X}(\Sigma)$ of the solution $u$.
In the following, we consider %
coefficients $\opA_j\in L^\infty(\R^\dim;\mathbb S^\dim)$ 
and $p_j\in L^\infty(\R^\dim)$ that satisfy \Cref{ass:A_p_new} and 
agree with $\opA$ and $p$ on $\Omega_j$
(but possibly differ on $\R^\dim\setminus\overline{\Omega_j}$), i.e.,
\begin{align}\label{eqn:opAj_pj_def}
	\opA_j|_{\Omega_j}=\opA|_{\Omega_j}
	\qquad\text{and}\qquad
	p_j|_{\Omega_j}=p|_{\Omega_j}
	\qquad\text{for }j=0,\dots,J.
\end{align}
The corresponding single layer, double layer, and Calder\'on operators from 
\Cref{sec:Potential operators for interface problems} (for $G$ replaced by $\Omega_j$)
are denoted by $\opSj$, $\opDj$, and $\opCj$.
(Recall for piecewise Lipschitz $\opA_j$ that \Cref{ass:A_p_new} reduces to%
~\ref{ass:C3}--\ref{ass:C4} 
by \Cref{lem:uniqueness}.)
\Cref{rem:Aj_pj} below discusses the freedom to define the coefficients $\opA_j$ and $p_j$ on $\Omega\setminus\overline\Omega_j$.
The Calder\'on identity in \Cref{lem:C_equiv}
relates 
the solution $u\in H^1_{\mathrm{loc}}(\Omega\setminus\Sigma)$
to~\eqref{eqn:ATP} with its Cauchy traces 
$\traceC u=\umult_{\Sigma}=(\umult_{\Sigma,j})_{j=0}^J\in\mathbb{X}(\Sigma)$
by (cf.\ \Cref{thm:TP_equivalence} for details)
\begin{align*}
	-\Div(\opA_j\nabla u) +s^2p_j u &= 0\quad\text{in }\Omega_j	
	\quad\Leftrightarrow\quad
	(\opCj-\tfrac12)\umult_{\Sigma,j} = 0,
	\qquad\text{for all }j=0,\dots,J.
\end{align*}
This leads to a multi-trace formulation of~\eqref{eqn:ATP} 
that seeks $\umult_{\Sigma}=(\umult_{\Sigma,j})_{j=0}^J\in\mathbb{X}(\Sigma)$ with
\begin{equation}\label{eqn:MTP}
\boxed{
	\begin{aligned}
		(\opCj -\tfrac12)\umult_{\Sigma,j} 
		&=0 \qquad\text{for all }j=0,\dots,J,\\
		\umult_{\Sigma}-\bbg&\in\mathbf{X}_0(\Sigma).
	\end{aligned}
}
\end{equation}
The substitution of 
$\bbt=\umult_{\Sigma}-\bbg\in \mathbf{X}_0(\Sigma)$ in~\eqref{eqn:MTP}
results in a single-trace formulation of the transmission problem~\eqref{eqn:ATP}
that seeks $\bbt\in \mathbf{X}_{0}(\Sigma)$ with
\begin{equation}\label{eqn:S_ATP}
\boxed{
	\begin{aligned}
		\left(\opCj-\tfrac{1}{2}\right)\bbt[\Sigma,j]
			&= -\left(\opCj-\tfrac{1}{2}\right)\bbg[\Sigma,j]
				&&\text{for } j=0,\dots,J.\\
	\end{aligned}
}
\end{equation}
Recall $\slow$ and $\sLow=\slow^{-1}$ from~\eqref{eqn:s_low_def} 
and that the full-space Helmholtz problem~\eqref{eqn:FSP}
for $s\in\mathbb{C}_{\geq0}^*$,
$\opA\in L^\infty(\Omega;\mathbb{S}^n), p\in L^\infty(\Omega)$ with~\ref{ass:C1}--\ref{ass:C3}
satisfies \Cref{ass:A_p}.
\begin{theorem}[equivalence for the transmission problem]\label{thm:TP_equivalence}
	For any $\bbg\in\mathbb{X}(\Sigma)$, the solutions
	$u\in H^{1}_{\mathrm{loc}}(\Omega\setminus\Sigma,\opA)$ to~\eqref{eqn:ATP},
	$\umult_{\Sigma}=(\umult_{D,j},\umult_{N,j})_{j=1}^J\in\Xmult(\Sigma)$
	to~\eqref{eqn:MTP}, and
	$\bbt\in \mathbf{X}(\Sigma)$ to~\eqref{eqn:S_ATP}
	exist uniquely and satisfy
	\begin{enumerate}[label=(\roman*)]
	\item $u|_{\Omega_j} = (\opSj \umult_{N,j}-\opDj\umult_{D,j})|_{\Omega_j}$ for all $j=0,\dots,J$,
	\item $\umult_{\Sigma} = \traceC u$ and $\bbt= \traceC u - \bbg$,
	\item $\displaystyle
		\Const{ap}^{-1}\slow\,\|\umult_{\Sigma}\|_{\mathbb{X}(\Sigma),s}
		\leq\|u\|_{H^1(\Omega\setminus\Sigma),s}
		\leq \Const{ap}(1+\CN(s))\sLow\,
		\inf_{\bbh[]\in\mathbf{X}_0(\Sigma)}\|\bbg-\bbh[]\|_{\mathbb{X}(\Sigma),s}$.
	\end{enumerate}
	The constant $\Const{ap}>0$ is independent of $s$ and exclusively depends 
	on $a_{\max}$, $p_{\max}$, and on $\Omega_j$ for
	$j=0,\dots,J$.
\end{theorem}
\begin{proof}
The proof of~\Cref{thm:TP_equivalence} splits into three steps.

\medskip
\noindent\emph{Step 1} deduces the equivalence of the three formulations 
	from \Cref{lem:C_equiv}, with similar arguments as 
	in~\cite{Von:BoundaryIntegralEquations1989} given here for completeness.

	\medskip
	\eqref{eqn:ATP}$\Rightarrow$\eqref{eqn:MTP} with~\emph{(ii)}:
	Let $u\in H^1_{\mathrm{loc}}(\Omega\setminus\Sigma)$ be a solution 
	to~\eqref{eqn:ATP} and
	define
	$u_j\in H^1(B_R\setminus\Gamma_j)$ by 
	$u_j\coloneqq u|_{\Omega_j\cap B_R}$ and
	$u_j|_{B_R\setminus\overline{\Omega_j}}\equiv 0$ for all $j=0,\dots,J$.
	By construction, $u_j\in V(B_R\setminus\Gamma_j,\opA_j,s)$ solves
	the interface problem
	\begin{equation}\label{eqn:u_j_interface}
	\begin{aligned}
		-\Div(\opA_j\nabla u_j) + s^2p_j\,u_j
			&=0&&\text{in }B_R\setminus\Gamma_j,\\
		\jumpD[\Omega_j]{u_j} = \traceD[\Omega_j]u
		\quad\text{and}\quad
		\jumpN[\Omega]{u_j} &= \traceN[\Omega_j]u&&\text{on }\Gamma_j.
	\end{aligned}
	\end{equation}
	Hence 
	\Cref{lem:C_equiv} for 
	$\umult_{\Sigma,j}\coloneqq(\traceD[\Omega_j] u,\traceN[\Omega_j] u)
	\in \mathbb{X}(\Gamma_j)$
	and $\traceDext[\Omega_j]u_j =0=\traceNext[\Omega_j]u_j$
	reveals $(\opCj -\tfrac12)\umult_{\Sigma,j}=0$.
	This and 
	$\traceC u = \umult_{\Sigma}\coloneqq (\umult_{\Sigma,j})_{j=1}^J\in\Xmult(\Sigma)$
	verify that $\umult_{\Sigma}$
	solves~\eqref{eqn:MTP}.

	\medskip
	\eqref{eqn:MTP}$\Rightarrow$\eqref{eqn:ATP} with \emph{(i)}.
	Let $\umult_{\Sigma}=(\umult_{D,j},\umult_{N,j})_{j=0}^J\in\Xmult(\Sigma)$ solve
	\eqref{eqn:MTP}, set
	\begin{align*}
		u_j\coloneqq \opSj \umult_{\mathrm{N},j} - \opDj\umult_{\mathrm{D},j}
		\in V(B_R\setminus\Gamma_j,\opA_j,s)
		\qquad\text{for all }j=0,\dots,J,
	\end{align*}
	and define $u\in H^1_{\mathrm{loc}}(\Omega\setminus\Sigma)$ by \emph{(i)}.
	\Cref{lem:C_equiv} and~\eqref{eqn:MTP} verify
	$\umult_{\Sigma,j}=(\traceD[\Omega_j]u_j,\traceN[\Omega_j]u_j)$.
	Since $u|_{\Omega_j}= u_j|_{\Omega_j}$ by construction, 
	this implies $\traceC u=\umult_{\Sigma}$.
	Green's representation formula in \Cref{lem:integral_formulation}
	implies that $u_j$ solves the interface problem~\eqref{eqn:u_j_interface}.
	Hence $u$ solves~\eqref{eqn:ATP}.

	\medskip
	The equivalence of the formulations
	\eqref{eqn:MTP}$\Leftrightarrow$\eqref{eqn:S_ATP} 
	with $\umult_{\Sigma}=\bbt+\bbg$ is obvious.

	\medskip
	\noindent\emph{Step 2} establishes the existence and uniqueness.
	Let $\bbg=(\bbg[\mathrm{D},j],\bbg[\mathrm{N},j])_{j=0}^J\in \mathbb{X}(\Sigma)$ 
	be arbitrary and consider the equivalent reformulation of~\eqref{eqn:ATP}
	(cf.\ \Cref{sec:potential_operators,sec:Potential operators for interface problems})
	on the truncated domain $\Omega_R=\Omega\cap B_R$
	that seeks $u\in H^1(\Omega_R)$ with
	\begin{equation}\label{eqn:TATP}
		\begin{aligned}
			-\Div(\opA \nabla u) +s^2 p\, u
				&= 0 
					&&\text{in }\Omega_R\setminus\Sigma,\\
			\partial_r u&=\DtN u&&\text{on }S_R,\\
			\traceC u-\bbg&\in\mathbf{X}_0(\Sigma).
		\end{aligned}
	\end{equation}
	By the surjectivity of the trace operator,
	there exists $u_{\Sigma}\in H^1(\Omega_R\setminus\Sigma)$ with 
	$\traceD[\Omega_j]u_{\Sigma} = \bbg[\mathrm{D},j]$ for $j=0,\dots,J$
	and we may and will assume $\mathrm{supp}u_{\Sigma}\subset B_R$ such that
	$\partial_r u_{\Sigma}=0=\DtN u_{\Sigma}$ on $S_R$ in the following.
	Recall the sesquilinear form $\opl(\bullet,\bullet)$ 
	associated to the Helmholtz operator
	on $\Omega_R$ with DtN boundary conditions on $S_R$ from~\eqref{eqn:l_def}.
	The weak formulation of~\eqref{eqn:TATP} seeks $u=u_0+u_{\Sigma}$ with
	$u_0\in H^1_{\Gamma_{\mathrm{D}}}(\Omega_R)$ and
	\begin{align}\label{eqn:W_TATP}
		\opl(u_0,\varphi) 
		=F(\conj \varphi)
	\end{align}
	for all $\varphi\in H^1_{\Gamma_{\mathrm{D}}}(\Omega_R)$,
	where the right-hand side $F\in \Hdual_{\Gamma_{\mathrm{D}}}(\Omega_R)$ is given by
	\begin{align*}
		F(\conj \varphi) 
		\coloneqq \sum^{J}_{j=0} \left\langle \bbg[\mathrm{N},j],\conj\varphi\right\rangle_{\Gamma_j}
		- \int_{\Omega_R\setminus\Sigma}
			\left(\opA\nabla u_{\Sigma}\cdot \nabla \conj \varphi 
			+s^2p\,u_{\Sigma}\conj \varphi\right)\d x.
	\end{align*}
	Since~\eqref{eqn:W_TATP} is of the form~\eqref{eqn:WTHP},
	\Cref{thm:L_N_well_def} provides the existence of a solution 
	$u_0\in H^1_{\Gamma_{\mathrm{D}}}(\Omega_R)$
	to~\eqref{eqn:W_TATP} and we set $u\coloneqq u_0+u_{\Sigma}$.
	A piecewise integration by parts 
	(over $\Omega_0\cap B_R$ and $\Omega_1,\dots,\Omega_J$) for
	the right-hand side in~\eqref{eqn:W_TATP} 
	and $\partial_r u_{\Sigma}=0=\DtN u_{\Sigma}$ reveal
	\begin{align}\label{eqn:W_TATP_div}
		\int_{\Omega_R\setminus\Sigma}
		-\Div(\opA\nabla u+s^2p\,u) \conj\varphi\d x
		+ \left\langle \jumpNext u, \conj\varphi\right\rangle_{S_R}
		=\sum^{J}_{j=0} 
		\left\langle \bbg[\mathrm{N},j]
		-\traceN[\Omega_j]u,\conj\varphi\right\rangle_{\Gamma_j}
	\end{align}
	with $\jumpNext\bullet$ from~\eqref{eqn:V_BR_def}
	for all $\varphi\in H^1_{\Gamma_{\mathrm{D}}}(\Omega_R)$.
	This implies $-\Div(\opA \nabla u) + s^2p\,u=0$ in $\Omega_R\setminus\Sigma$
	and $\jumpNext u=0$ on $S_R$ 
	(by arguments similar to the proof of \Cref{lem:Newton_potential_L2}).
	Hence it remains to verify $\traceC u-\bbg\in \mathbf{X}_0(\Sigma)$.
	The construction of $u_{\Sigma}$ and $u=u_0-u_\Sigma|_{\Omega_R}$ show
	\begin{align*}
		\traceD[\Omega_j] u -\bbg[\mathrm{D},j]
		= \traceD[\Omega_j] u_0 + \traceD[\Omega_j] u_{\Sigma} - \bbg[\mathrm{D},j]
		= \traceD[\Omega_j] u_0
		\qquad\text{for all }j=0,\dots,J.
	\end{align*}
	In other words, the Dirichlet part of $\traceC u-\bbg$ 
	is the trace of some extension of $u_0$ to $H^1_{\Gamma_{\mathrm{D}}}(\Omega)$.
	Let $\mathbf{q}\in H(\Omega\setminus\Sigma,\Div)$ satisfy $\traceNu[\Omega_j]\mathbf{q} =
	\traceN[\Omega_j]u-\bbg[\mathrm{N},j]$ for all $j=0,\dots,J$ and
	$\mathrm{supp}(\mathbf{q})\subset B_R$.
	Since $u$ solves the Helmholtz equation in $\Omega_R\setminus\Sigma$ with 
	$\jumpNext u =0$ on $S_R$, the left-hand side in~\eqref{eqn:W_TATP_div} vanishes.
	This and a piecewise integration by parts
	result in
	\begin{align*}
		0 = \sum^{J}_{j=0} 
		\left\langle \traceN[\Omega_j]u-\bbg[\mathrm{N},j],
			\conj\varphi\right\rangle_{\Gamma_j}
			= \sum^{J}_{j=0} \left\langle \traceNu[\Omega_j] \mathbf{q}, 
			\conj \varphi\right\rangle_{\Gamma_j}
			=\int_{\Omega_R\setminus\Sigma}\left(\mathbf{q}\cdot\nabla\conj\varphi 
				+ \conj \varphi\Div \mathbf{q}\right)\d x
	\end{align*}
	for all $\varphi\in H^1_{\Gamma_{\mathrm{D}}}(\Omega_R)$.
	Hence (by definition of the weak divergence)
	$\Div \mathbf{q}\in L^2(\Omega_R)$ and $\mathbf{q}\cdot\nu_{\Gamma}=0$ on $\Gamma_{\mathrm{N}}$.
	This and $\mathrm{supp} (\mathbf{q})\subset B_R$
	imply $\mathbf{q}\in H_{\Gamma_{\mathrm{N}}}(\Omega,\Div)$.
	Consequently, $\traceC u-\bbg\in \mathbf{X}_0(\Sigma)$.

	Since the solution $u_0\in H^1_{\Gamma_{\mathrm{D}}}(\Omega_R)$ 
	to~\eqref{eqn:W_TATP} is in fact unique by
	\Cref{thm:L_N_well_def}, the previous arguments verify the existence of 
	a unique solution to~\eqref{eqn:TATP} (and equivalently~\eqref{eqn:ATP}).
	With the already establised equivalence of~\eqref{eqn:ATP}
	to the multi-trace and single-trace formulations~\eqref{eqn:MTP}--\eqref{eqn:S_ATP}
	with \Cref{thm:TP_equivalence}.i--iii, this concludes Step 2.

	\medskip
	\noindent\emph{Step 3} provides wavenumber-explicit bounds 
	based on a particular Dirichlet lifting $u_{\Sigma}$ as used in Step 2.
	The properties of
	minimal extension norms (cf.~\cite[Thm.~3.1]{Gra:OptimalTraceNorms2025}) 
	imply that the minimum in~\eqref{eqn:H12_norm_def} is attained. 
	Hence there is $v_{\Sigma}\in H^1(\Omega\setminus\Sigma)$ with 
	\begin{align*}
		\traceD[\Omega_j]{v_{\Sigma}} = \bbg[\mathrm{D},j]
		\qquad\text{and}\qquad
		\NormHD[\Gamma_j]{s}{\bbg[\mathrm{D},j]} = \|v_{\Sigma}\|_{H^1(\Omega_j),s}
	\end{align*}
	for all $j=0,\dots,J$.
	Set $u_{\Sigma} \coloneqq \varphi v_{\Sigma}\in H^1(\Omega\setminus\Sigma)$ 
	for some $\varphi\in C^\infty_0(B_R)$ with $\varphi\equiv 1$ on $\Sigma$.
	Observe that $u_{\Sigma}$ satisfies the properties from Step 2
	and, by the product rule,
	\begin{align}\label{eqn:gD_bound}
		\Const{\varphi}^{-1}\|u_{\Sigma}\|_{H^1(\Omega_R\setminus\Sigma),s}
		\leq
		\|v_{\Sigma}\|_{H^1(\Omega\setminus\Sigma),\max\{1,s\}}
		\leq\sLow\sqrt{\sum^{J}_{j=0} \NormHD[\Gamma_j]{s}{\bbg[\mathrm{D},j]}^2}
	\end{align}
	for $\Const{\varphi}
	\coloneqq\|\varphi\|_{L^\infty(\Omega_R)}+\|\nabla\varphi\|_{L^\infty(\Omega_R)}$.
	By Step 2, the unique solution $u\in H^1(\Omega\setminus\Sigma)$ to~\eqref{eqn:ATP}
	splits as $u= u_0+u_{\Sigma}$ with the 
	solution $u_0\in H^1_{\Gamma_{\mathrm{D}}}(\Omega_R)$ to~\eqref{eqn:W_TATP}.
	Cauchy inequalities and \eqref{eqn:Ctr_D} in \Cref{lem:trace_inequality}
	control the dual norm~\eqref{eqn:H1_dualnorm_def}
	of the right-hand side in~\eqref{eqn:W_TATP} by
	\begin{align*}
		\|F\|_{\Hdual_{\Gamma_{\mathrm{D}}}(\Omega_R)}
		\leq
		\sqrt{\sum^{J}_{j=0} \NormHN[\Gamma_j]{s}{\bbg[\mathrm{N},j]}^2}
			+\max\{a_{\max},p_{\max}\}\|u_{\Sigma}\|_{H^1(\Omega_R\setminus\Sigma),s}
	\end{align*}
	This, \Cref{thm:L_N_well_def}, and~\eqref{eqn:gD_bound} combined with 
	triangle and Cauchy inequalities provide
	\begin{align}\nonumber
		\|u\|_{H^1(\Omega_R\setminus\Sigma),s}
		&\leq
		\CN(s) \|F\|_{\Hdual_{\Gamma_{\mathrm{D}}}(\Omega_R)}
		+ \|u_{\Sigma}\|_{H^1(\Omega\setminus\Sigma),s}\\
		&\leq\Const{ap,1}(1+\CN(s))\sLow\,\|\bbg\|_{\mathbb{X}(\Sigma),s}
		\label{eqn:u_norm_bound}
	\end{align}
	for a constant $\Const{ap,1}$ that exclusively depends on $\Omega_R\setminus\Sigma$ 
	(through $\Const{\varphi}$), $a_{\max}$ and $p_{\max}$.

	By \Cref{lem:trace_inequality}, the Cauchy trace operator $\traceC$ satisfies
	\begin{align*}
		\|\traceD[\Omega_j] u\|_{H^{1/2}(\Gamma_j),s}
		&\leq \|u\|_{H^1(\Omega_j\cap B_R),s}\\
		\NormHN[\Gamma_j]{s}{\traceN[\Omega_j]{u}}
		&\leq \Const{tr,N}\sLow\, \|\opA\nabla u\|_{H(\Omega_j\cap B_R,\Div),s},
	\end{align*}
	where $\Const{tr,N}\sLow$ may be replaced by $1$ for $j=1,\dots,J$, but not for 
	$j=0$.
	Using $\Div(\opA \nabla u) = s^2 p u$ in $\Omega_R\setminus\Sigma$
	to control
	$\|\opA\nabla u\|_{H(\Omega_R\setminus\Sigma,\Div),s}$ by
	$\|u\|_{H^1(\Omega_R\setminus\Sigma),s}$, this implies
	\begin{align}\label{eqn:traceC_bound}
		\|\traceC u\|_{\mathbb{X}(\Sigma),s}
		\leq \Const{ap,2}\sLow\,\|u\|_{H^1(\Omega\setminus\Sigma),s}
	\end{align}
	for some constant $\Const{ap,2}>0$
	that exclusively depends on $\Const{tr,N}$, $a_{\max}$, and $p_{\max}$.
	Since the solution $u$ to~\eqref{eqn:ATP} 
	(and $\umult_{\Sigma}=\traceC u$ by \emph{(ii)}) 
	is the same for all transmission data in 
	$\{\bbg+\bbh[]\ :\ \bbh[]\in \mathbf{X}_0(\Sigma)\}$, the combination%
	~\eqref{eqn:u_norm_bound}--\eqref{eqn:traceC_bound} results in \emph{(iii)}
	for $\Const{ap}\coloneqq\max\{\Const{ap,1},\Const{ap,2}\}$
	and concludes the proof.
\end{proof}

\begin{remark}\label{rem:Aj_pj}
	The choice $\opA_j=\opA$ and $p_j=p$ is allowed in~\eqref{eqn:opAj_pj_def}
	but \emph{not} required to define $\opSj, \opDj$, and $\opCj$ 
	for the subsets $\Omega_j$ for $j=0,\dots,J$. 
	\Cref{thm:TP_equivalence}.i reveals that $\opSj$ and $\opDj$
	only affect the solution on $\Omega_j$ where $\opA_j=\opA$ and $p_j=p$ hold by~\eqref{eqn:opAj_pj_def}.

For a piecewise constant, isotropic coefficient $\mathbb{A}$
and piecewise constant $p$ in~\eqref{HelmEq1}, a typical choice of the
subdomains $\Omega_{j}$ ensures $\left.  \mathbb{A}\right\vert
_{\Omega_{j}}=c_{j}\mathbb{I}$ and $\left.  p\right\vert _{\Omega_{j}}=p_{j}$
for some positive constants $c_{j},p_{j}\in\mathbb{R}$.
Without loss of generality, 
we may assume $c_{0}=1=p_{0}$ by a simple scaling of the subproblems on $\Omega_j$,
so that~\ref{ass:C3} is satisfied. 
The constant extensions $\mathbb{A}_{j}$ and $p_{j}$ of 
$\left.  \mathbb{A}\right\vert _{\Omega_{j}}$ and $\left.
p\right\vert _{\Omega_{j}}$ to $\mathbb{R}^{n}$ lead to \textit{globally}
constant coefficients $\mathbb{A}_{j}$ and $p_{j}$, which in general
differ from the original \textit{piecewise} constant coefficients
$\mathbb{A}$ and $p$ in the transmission problem~\eqref{eqn:ATP}.
In this setting ($\mathbb{A}_{j}$ and $p_{j}$ globally constant), 
the single and double layer operator as well as the associated
skeleton operators admit classical representations as boundary integral
operators with explicitly known kernel functions%
~\cite{McL:StronglyEllipticSystems2000,SS:BoundaryElementMethods2011}. 
In general, the flexibility in the choice of the extensions $\mathbb{A}_{j}$, $p_{j}$ 
can be used advantageously for certain piecewise (Lipschitz) smooth coefficients 
to derive representation as boundary integral operators, 
e.g., by asymptotic methods%
~\cite{BDM:Asymptotics2DWhispering2021,BT:AsymptoticGreenFunction2010,JT:MatchingAsymptoticExpansions2006}, 
Luneburg spheres~\cite{LLA:ScatteringPlaneElectromagnetic2015},
or WKB methods~\cite{Goo:RoleAttractivePotential1971,Eng:ApproximateGreenFunction1983}.
We emphasise that, regardless of the representation as integral operators,
the skeleton integral method in this paper provides a stable variational formulation 
for these non-local operator equations.
\end{remark}

\begin{remark}
	The transmission problem~\eqref{eqn:ATP}
	can be interpreted as a special case of the exterior problem~\eqref{eqn:FSP}
	for the non-constant (e.g., piecewise Lipschitz) coefficient matrix $\opA\in L^\infty(\Omega,\mathbb{S}^\dim)$ (and
	$p\in L^\infty(\Omega)$).
	Indeed, the proof of \Cref{thm:TP_equivalence} reveals that the equivalent weak formulation~\eqref{eqn:W_TATP}
	of the transmission problem is of the form~\eqref{eqn:WTHP} for some right-hand side $F\in
	\Hdual_{\Gamma_\mathrm{D}}(\Omega_R)$ that represents the transmission data $\bbg\in\mathbb{X}(\Sigma)$.
\end{remark}
\subsection{Variational single-trace formulation}%
\label{sub:Variational formulation}

The single-trace formulation~\eqref{eqn:S_ATP} can be rewritten as the operator equation
\begin{align*}
	\opAvar \bbt = -\opAvar\bbg
\end{align*}
for the given data $\bbg\in\mathbb{X}(\Sigma)$, the solution $\bbt\in\mathbf{X}_0(\Sigma)$,
and the operator
\begin{align}\label{eqn:Avar_dev}
	\opAvar\coloneqq -\tfrac12 \operatorname{id}+ \mathrm{diag}(\opCj\ :\ j=0,\dots,J)
	:\mathbb{X}(\Sigma)\to \mathbb{X}(\Sigma)
	.
\end{align}

\Cref{thm:TP_equivalence} states that $\opAvar$ is an isomorphism from $\mathbf{X}_0(\Sigma)$ 
onto the (total/full) image $\opAvar \mathbb{X}(\Sigma)$.
In partiular, $\opAvar\mathbb{X}(\Sigma)=\opAvar\mathbf{X}_0(\Sigma)$ and%
~\eqref{eqn:S_ATP} may be recasted in a variational least-squares formulation:
Given $\bbg\in\mathbb{X}(\Sigma)$, find $\bbt\in\mathbf{X}_0(\Sigma)$ with
\begin{align}\label{eqn:STF_var}
	\boxed{
	\left\langle \opAvar\bbt,\conj{\opAvar\bbh}\right\rangle_{\mathbb{X}(\Sigma)}
	=-\left\langle \opAvar\bbg,\conj{\opAvar\bbh}\right\rangle_{\mathbb{X}(\Sigma)}
	\qquad\text{for all }\bbh\in\mathbf{X}_0(\Sigma).
}
\end{align}
Note, that the operator $\opAvar$ acts on the test \emph{and} trial functions in~\eqref{eqn:STF_var}.
This is not ideal from a practical (implementation) point of view.
Motivated by the isomorphism between $\mathbf{X}_0(\Sigma)$ and $\opAvar\mathbb{X}(\Sigma)$,
we consider the following alternative problem
in the spirit of the classical single-trace formulation of first
kind~\cite[Sec.~4]{CH:IntegralEquationsAcoustic2015}.
Given $\bbg\in\mathbb{X}(\Sigma)$, seek a solution $\bbt\in\mathbf{X}_0(\Sigma)$ to%
\begin{align}\label{eqn:STF_classical_var}
	\boxed{
	\left\langle \opAvar\bbt,\conj{\bbh}\right\rangle_{\mathbb{X}(\Sigma)}
	=-\left\langle \opAvar\bbg,\conj{\bbh}\right\rangle_{\mathbb{X}(\Sigma)}
	\qquad\text{for all }\bbh\in\mathbf{X}_0(\Sigma).
}
\end{align}%
\begin{theorem}[variational single-trace formulations]\label{thm:STP_var}
	For all $\bbg\in\mathbb{X}(\Sigma)$,
	the unique solution $\bbt\in\mathbf{X}_0(\Sigma)$ to~\eqref{eqn:S_ATP}
	\begin{itemize}
		\item is the unique solution to~\eqref{eqn:STF_var},
		\item solves~\eqref{eqn:STF_classical_var} and the solution set 
			of~\eqref{eqn:STF_classical_var} is given by 
			$\bbt+\mathrm{ker}_{\mathbf{X}_0(\Sigma)}\opAvar$ with
			\begin{align*}
				\mathrm{ker}_{\mathbf{X}_0(\Sigma)}\opAvar
				\coloneqq \left\{\bbg[]\in\mathbf{X}_0(\Sigma)\ :\ \left\langle
					\opAvar\bbg[],\conj{\bbh[]}\right\rangle_{\mathbb{X}(\Sigma)}=0
				\quad\text{for all }\bbh[]\in \mathbf{X}_0(\Sigma)\right\}.
			\end{align*}
	\end{itemize}
\end{theorem}
\begin{proof}
	Let $\bbt\in\mathbf{X}_0(\Sigma)$ denote the unique solution to~\eqref{eqn:S_ATP}, i.e., 
	$\opAvar\bbt=-\opAvar\bbg$. It is clear that $\bbt$ solves~\eqref{eqn:STF_var}
	as well as~\eqref{eqn:STF_classical_var}.
	For any other solution $\widetilde{\bbt}\in\mathbf{X}_0(\Sigma)$ to~\eqref{eqn:STF_var},
	$\bbh\coloneqq\bbt-\widetilde{\bbt}\in\mathbf{X}_0(\Sigma)$ is an admissible test function
	for~\eqref{eqn:STF_var}. Hence $\|\opAvar(\bbt-\widetilde{\bbt})\|_{\mathbb{X}(\Sigma),s}=0$
	and \Cref{thm:TP_equivalence} implies $\bbt=\widetilde{\bbt}$.
	The second statement is clear.
\end{proof}
\Cref{thm:STP_var} states that the well-posedness of~\eqref{eqn:STF_var} is unconditionally equivalent
to~\eqref{eqn:S_ATP}, 
while~\eqref{eqn:STF_classical_var} is equivalent to~\eqref{eqn:S_ATP} if and only if 
$\mathrm{ker}_{\mathbf{X}_0(\Sigma)}\opAvar=\{0\}$ is trivial.
The latter always holds if the wavenumber has positive real part $\Re s>0$.
Indeed, the operator $\opAvar$ (for $\Re s>0$) is coercive
\begin{align*}
	\Re\left\langle \opAvar\bbg[],\conj{\bbg[]}\right\rangle 
	\geq c_{\mathrm{A}}(s)\|\bbg[]\|_{\mathbb{X}(\Sigma),s}^2
	\qquad\text{for all }\bbg[]\in\mathbb{X}(\Sigma).
\end{align*}
(This follows with minor modifications as 
\cite[Thm.~34]{FHS:SkeletonIntegralEquations2024}, where the case $R=\infty$ is discussed.)
\Cref{rem:kernel_degeneracy} below discusses the case of purely imaginary wavenumbers 
($\Re s=0$), where 
$\mathrm{ker}_{\mathbf{X}_0(\Sigma)}\opAvar\ne\{0\}$ is possible and related to 
certain geometrical configurations of the scatterer and the interfaces.
In any case, the dimension of $\mathrm{ker}_{\mathbf{X}_0(\Sigma)}\opAvar$ is finite
as the interpretation 
of the following G\r arding-type inequality in \Cref{rem:STF_var} below shows.
\begin{lemma}[G\r arding inequality]\label{lem:STF_var_Garding}
	The operator $\opAvar:\mathbb{X}(\Sigma)\to\mathbb{X}(\Sigma)$ is continuous.
	There exists a compact operator 
	$\mathsf T(s):\mathbf{X}_0(\Sigma)\to \mathbb{X}(\Sigma)$ such that 
	\begin{align*}
		\Re\left\langle (\opAvar+\mathsf T(s)) \bbg[] ,
		\conj{\bbg[]} \right\rangle_{\mathbb{X}(\Sigma)}
		&\geq c_{\mathrm{A}} \,\slow^2\,
		\|\bbg[] \|_{\mathbb{X}(\Sigma),s}^2
		\qquad\text{for all }\bbg\in\mathbf{X}_0(\Sigma).
	\end{align*}
	The constant $c_{\mathrm{A}}>0$ exclusively depends on $\Sigma,R,\opAvar$, and $p$.
\end{lemma}
\begin{proof}
	The continuity of $\opAvar$ is clear from that of the Calder\'on operators (by \Cref{lem:C_garding}).
	Let $\mathsf T_j(s):\mathbf{X}(\Gamma_j)\to \mathbf{X}(\Gamma_j)$ denote the 
	compact operator from \Cref{lem:C_garding}
	for all $j=0,\dots,J$.
	It is well-known (e.g., by%
	~\cite[Lem.~4.1]{CH:IntegralEquationsAcoustic2015} and%
	~\cite[Rem.~5.5]{CHJP:NovelMultitraceBoundary2015})
	that the self-dual pairing~\eqref{eqn:X_dual} satisfies
	\begin{align*}
		0=\sum^{J}_{j=0} \left\langle \bbg[j],\conj{\bbh[j]}\right\rangle_{\mathbf{X}(\Gamma_j)} 
		\qquad\text{for all }\bbg[j],\bbh[j]\in\mathbf{X}_0(\Sigma).
	\end{align*}
	This and the coercivity by \Cref{lem:C_garding} reveal for
	$\mathsf T(s)\coloneqq \mathrm{diag}(\mathsf T_j(s)\ :\ j=0,\dots,J)$ that
	\begin{align*}
		\Re\left\langle (\opAvar+\mathsf T(s))\bbg[],\conj{\bbg[]}\right\rangle_{\mathbb{X}(\Sigma)}
		&= \sum^{J}_{j=0} 
		\Re\left\langle (\opCj+\mathsf T_j(s))\bbg[j],\conj{\bbg[j]}\right\rangle_{\mathbf{X}(\Gamma_j)}\\
		&\geq
		c_{\mathrm{A}} \,\slow^2\,
		\sum^{J}_{j=0}\|\bbg[j] \|_{\mathbf{X}(\Gamma_j),s}^2
	\end{align*}
	for some constant $c_{\mathrm{A}}>0$ that exclusively depends on $R$ and $\opA_j, \Gamma_j$ for $j=0,\dots,J$.
\end{proof}

\begin{remark}[operator equations]\label{rem:STF_var}
	The variational formulations~\eqref{eqn:STF_var} with 
	$R(\Sigma)\coloneqq \opAvar\mathbb{X}(\Sigma)$ and~\eqref{eqn:STF_classical_var}
	with $R(\Sigma)\coloneqq\mathbf{X}_0(\Sigma)$ can be written in operator notation as
	\begin{align}\label{eqn:STF_var_operator}
		P_{R(\Sigma)}\opAvar \bbt = - P_{R(\Sigma)}\opAvar\bbg
	\end{align}
	using the projection $P_{R(\Sigma)}:\mathbb{X}(\Sigma)\to R(\Sigma)\subset\mathbb{X}(\Sigma)$ 
	defined by
	\begin{align*}
		\left\langle P_{R(\Sigma)}\bbg[],\conj{\bbh[]}\right\rangle_{\mathbb{X}(\Sigma)}
		= \left\langle \bbg[],\conj{\bbh[]}\right\rangle_{\mathbb{X}(\Sigma)}
		\qquad\text{for all }\bbg[]\in\mathbb{X}(\Sigma), \bbh[]\in R(\Sigma).
	\end{align*}
	Since $P_{\opAvar\mathbb{X}(\Sigma)}$ is the identity on the image of $\opAvar$,%
	~\eqref{eqn:STF_var_operator} coincides with~\eqref{eqn:S_ATP} for 
	$R(\Sigma)=\opAvar\mathbb{X}(\Sigma)$.
	For $R(\Sigma)=\mathbf{X}_0(\Sigma)$,%
	~\eqref{eqn:STF_var_operator} and~\eqref{eqn:S_ATP}
	are equivalent if and only if
	$P_{\mathbf{X}_0(\Sigma)}$ is invertible on the image of $\opAvar$.
	This is the case if
	\begin{align}\label{eqn:kernel_identity}
		\mathrm{ker}\,P_{\mathbf{X}_0(\Sigma)}\opAvar|_{\mathbf{X}_0(\Sigma)}
		=\mathrm{ker}_{\mathbf{X}_0(\Sigma)}\,\opAvar
		=\{0\}.
	\end{align}
	\Cref{lem:STF_var_Garding} implies that the operator
	$P_{\mathbf{X}_0(\Sigma)}(\opAvar + \mathsf T(s))|_{\mathbf{X}_0(\Sigma)}:
	\mathbf{X}_0(\Sigma)\to\mathbf{X}_0(\Sigma)$ is an isomorphism.
	Since $\mathsf T(s)$ is compact, 
	this means that $P_{\mathbf{X}_0(\Sigma)}\opAvar|_{\mathbf{X}_0(\Sigma)}$ is a
	Fredholm operator of index $0$ and has, in particular,
	a finite dimensional kernel%
	~\cite{Zei:NonlinearFunctionalAnalysis1992,Yos:FunctionalAnalysis1995,McL:StronglyEllipticSystems2000}.
\end{remark}
\begin{remark}[degeneracy of the kernel]\label{rem:kernel_degeneracy}
	In the critical case $\Re s=0$, 
	the kernel of $P_{\mathbf{X}_0(\Sigma)}\opAvar|_{\mathbf{X}_0(\Sigma)}$
	(i.e., $\mathrm{ker}_{\mathbf{X}_0(\Sigma)}\,\opAvar$ by~\eqref{eqn:kernel_identity})
	may be non-trivial and, by \Cref{thm:STP_var}, 
	the solutions to~\eqref{eqn:S_ATP} non-unique.
	This has been observed in%
	~\cite[Sec.~4]{CH:IntegralEquationsAcoustic2015} for piecewise constant coefficients
	and related to certain geometrical configurations of the scatterer 
	and the interfaces.
	The arguments therein 
	apply in the present case and reveal the analogue to%
	~\cite[Thm.~4.8]{CH:IntegralEquationsAcoustic2015}:
	$\mathrm{ker}_{\mathbf{X}_0(\Sigma)}\,\opAvar$ is non-trivial if and only if
	$\partial\Omega\subset \partial\Omega_j$ is a boundary component of some subdomain 
	$\Omega_j$ for $j=0,\dots,J$ and $\kappa_j=s p_j$ is an eigenvalue of the Laplacian 
	on the acoustic obstacle $\R^n\setminus\overline{\Omega}$ with
	Dirichlet and Neumann conditions on $\Gamma_{\mathrm{D}}$ and $\Gamma_{\mathrm{N}}$.
\end{remark}

\newpage
	\bibliographystyle{alphaabbr}
	\bibliography{Bibliography}

\end{document}